# Quantum categories, star autonomy, and quantum groupoids

Brian Day and Ross Street

**Abstract** A useful general concept of bialgebroid seems to be resolving itself in recent publications; we give a treatment in terms of modules and enriched categories. We define the term "quantum category". The definition of antipode for a bialgebroid is less resolved in the literature. Our suggestion is that the kind of dualization occurring in Barr's star-autonomous categories is more suitable than autonomy (= compactness = rigidity). This leads to our definition of quantum groupoid intended as a "Hopf algebra with several objects".

## 1. Introduction

This paper has several purposes. We wish to introduce the concept of quantum category. We also wish to generalize the theory of $*$-autonomous categories in the sense of [Ba1]. The connection between these two concepts is that they lead to our notion of quantum groupoid. It was shown by [Se] that $*$-autonomous categories provide models of the linear logic described in [Gi]. This suggests an interesting possibility of interactions between computer science and quantum group theory. Perhaps it will be possible, in future papers, to exploit the dichotomy between categories *as* structures and categories *of* structures. For example, what is the quantum category of finite sets, or the quantum category of finite dimensional vector spaces?

The section headings are as follows:

Before looking at quantum categories we will develop, in this introduction, a definition of "category" which suggests the definition of "quantum category". We will then relate this definition to the literature.

We use the terminology of Eilenberg-Kelly [EK] for monoidal categories and monoidal functors; so we use the adjective "*strong* monoidal" for a functor which preserves tensor and unit up to coherent natural isomorphisms. A comonoidal category would have, instead of a tensor product, a tensor coproduct $\mathcal{A} \longrightarrow \mathcal{A} \times \mathcal{A}$ and a counit with appropriately coherent constraints; this concept is not so interesting for ordinary categories but becomes more so for enriched categories. Comonoidal functors would go between comonoidal categories. So, for monoidal categories $\mathcal{A}$ and $\mathcal{X}$, like [McC], we use the term *opmonoidal functor* for a functor

$$F : \mathcal{A} \longrightarrow \mathcal{X}$$

equipped with natural family of morphisms $\delta_{A,B} : F(A \otimes B) \longrightarrow FA \otimes FB$ and a morphism $\varepsilon : FI \longrightarrow I$ that are coherent.

For any set $X$, consider the monoidal category $Set / X \times X$ of sets over $X \times X$ with the tensor product defined by

$$\left(A \xrightarrow{(s,t)} X \times X\right) \otimes \left(B \xrightarrow{(u,v)} X \times X\right) = \left(P \xrightarrow{(s \circ p, v \circ q)} X \times X\right)$$

where $P$ is the pullback of $t : A \longrightarrow X$ and $u : B \longrightarrow X$ with projections $p : P \longrightarrow A$ and $q : P \longrightarrow B$. The objects of $Set / X \times X$ are directed graphs with vertex-set $X$ and the monoids are the categories with object-set $X$; this is well known (see [ML]) and easy. Less well known, but also easy, is the fact that category structures on the graph $A \xrightarrow{(s,t)} X \times X$ amount to monoidal structures on the category $Set / A$ of sets over $A$ together with a strong monoidal structure on the functor

$$\Sigma_{(s,t)} : Set / A \longrightarrow Set / X \times X$$

defined on objects by composing the function into $A$ with $(s,t)$.

To see this, notice that every object of a slice category $Set / C$ is a coproduct of elements $c : 1 \longrightarrow C$ of $C$ (here $1$ is a chosen set with precisely one element) so that any tensor product on $Set / C$, which preserves coproducts in each variable, will be determined by its value on elements (which may not be another element in general). The tensor product on $Set / X \times X$ is such, and its value on elements is given by $(x,y) \otimes (u,v) = (x,v)$ when $y = u$ (which is in fact another element) but is the unique function $\varnothing \longrightarrow X \times X$ when $y \neq u$. Since $\Sigma_{(s,t)}$ is coproduct preserving and is to be strong monoidal, the tensor product on $Set / A$ preserves coproducts in each variable. An object of $Set / A$ has the same source set as its value under $\Sigma_{(s,t)}$. So, for elements $a$ and $b$ of $A$, the tensor product $a \otimes b$ is an element of $A$ if and only if $t(a) = s(b)$; in this case, $s(a \otimes b) = s(a)$ and $t(a \otimes b) = t(b)$; otherwise, $a \otimes b$ is the unique function $\varnothing \longrightarrow A$. The unit for the monoidal category $\mathrm{Set} / X \times X$ is the diagonal $X \longrightarrow X \times X$, so the unit for $\mathrm{Set} / A$ has the form $i : X \longrightarrow A$ with $s(i(x)) = t(i(x) = x$ for all $x \in X$. Thus we have a reflexive graph X, A, s, t, i with a composition operation. We leave the reader to check that the associativity and unity constraints of the monoidal category $\mathrm{Set} / A$ give associativity for the composition and that each $i(x)$ is an identity.

Conversely, suppose we have a category $\mathbf{A}$ with underlying graph $A \xrightarrow{(s,t)} X \times X$. Notice that $\mathrm{Set} / A$ is canonically equivalent to the category $[A, \mathrm{Set}]$ of functors from the discrete category A to Set. We can define a promonoidal structure (in the sense of [Da1]) on A by

$$P(a,b;c) = \begin{cases} 1 & \text{when } c \text{ is the composite of } a \text{ and } b; \\ \varnothing & \text{otherwise;} \end{cases} \qquad \text{and}$$

$$J(a) = \begin{cases} 1 & \text{when } a \text{ is an identity;} \\ \varnothing & \text{otherwise.} \end{cases}$$

Then $[A, \mathrm{Set}]$ becomes a monoidal category under convolution; this transports to a monoidal structure on $\mathrm{Set} / A$ for which $\Sigma_{(s,t)}$ is strong monoidal.

When our category $\mathbf{A}$ is actually a groupoid (that is, every arrow is invertible), there is a bijection $S : A \longrightarrow A$ defined by $Sa = a^{-1}$. We draw attention to the isomorphisms

$$P(a, b; Sc) \cong P(b, c; Sa),$$

noting here that S is its own inverse and that the diagram

$$\begin{array}{ccc} A & \xrightarrow{S} & A \\ {\scriptstyle (s,t)}\downarrow & & \downarrow{\scriptstyle (s,t)} \\ X \times X & \xrightarrow[S]{} & X \times X \end{array}$$

commutes, where the lower S is the switch map (which is inversion for X as a chaotic — or indiscrete — category). We will relate this kind of "antipode" structure to $*$-autonomy.

Now suppose we have a category $\mathbf{A} : A \xrightarrow{(s,t)} X \times X$ and suppose we regard $\mathrm{Set} / A$ as monoidal in the manner described above. The functor $\Sigma_{(s,t)}$ has a right adjoint $(s,t)^*$ defined by pulling back along $(s,t)$. The strong monoidal structure on $\Sigma_{(s,t)}$ is obviously both monoidal and opmonoidal; the opmonoidal structure transforms to a monoidal structure on the right adjoint $(s,t)^*$ in such a way that the unit and counit for the adjunction are monoidal natural transformations. The composite of monoidal functors is monoidal; so the endofunctor $G_{\mathbf{A}} = \Sigma_{(s,t)}(s,t)^*$ is also monoidal. The adjunction also generates a comonad structure on $G_{\mathbf{A}}$ in such a way that the counit and comultiplication are monoidal natural transformations; we have a monoidal comonad $G_{\mathbf{A}}$ on $\mathrm{Set}/X\times X$. Remember the term "monoidal comonad"!

It is also important to notice that $(s,t)^*$ has a right adjoint $\Pi_{(s,t)}$; so the endofunctor $G_{\mathbf{A}}$ has a right adjoint $(s,t)^*\ \Pi_{(s,t)}$. By Beck's Theorem (see [ML] for example), $\Sigma_{(s,t)}$ is comonadic since it is obviously conservative (that is, reflects isomorphisms) and preserves equalizers. On the other hand, any monoidal comonad on a monoidal category leads to a monoidal structure on the category of Eilenberg-Moore coalgebras in such a way that the forgetful functor is strong monoidal (see [McC] for example). Any cocontinuous endofunctor of $\mathrm{Set}/X\times X$ has the form $\Sigma_{(s,t)}(s,t)^*$ for some graph $A \xrightarrow{(s,t)} X\times X$. Assembling all this, we obtain:

**Proposition 1.1** *Categories with underlying graph* $A \xrightarrow{(s,t)} X\times X$ *are in bijection with monoidal comonad structures on the endofunctor* $\Sigma_{(s,t)}(s,t)^*$ *of* $\mathrm{Set}/X\times X$.

Let us compare the combinatorial context of Proposition 1.1 with the linear algebra context. Szlachányi [Szl] has shown that, for a k–algebra R, the $\times_R$-bialgebras of Takeuchi [Tak] are opmonoidal monads on the monoidal category $\mathrm{Vect}_k^{R\otimes R^\circ}$ of left R-, right R-bimodules over R where the underlying endofunctor of the monad is a left adjoint. These $\times_R$-bialgebras of Takeuchi have been convincingly proposed (see [Xu], [Lu], [Sch]) as the good concept of "bialgebroid" based on R (that is, with "object of objects R").

Here we face the usual dilemma. Given a k–bialgebra H, is it better to consider the category of modules for the underlying algebra with the monoidal structure coming from the comultiplication, or, the category of comodules for the underlying coalgebra with the monoidal structure coming from the multiplication? Our preference is definitely the latter since the obvious linearization of the group case leads to this decision; also see [JS2]. When H is finite dimensional (as a vector space over a field k) there is essentially no difference. Brzezinski-Militaru [BM] have already made the appropriate dualization of the $\times_R$-bialgebras of Takeuchi based on a k–coalgebra C rather than a k–algebra R. We take this as our concept of quantum category; it involves a monoidal comonad.

The basic examples of quantum groups are Hopf algebras with braidings (also called quasitriangular elements or R-matrices) or cobraidings, depending how the dilemma is resolved. Indeed, these basic quantum groups are cotortile bialgebras (see [JS2]). We leave it to a future paper to define and discuss braidings and twists on quantum categories.

So what is a quantum groupoid? It should be a quantum category with an "antipode". We first develop a notion of antipode for the $\times_R$-bialgebras of Takeuchi. We are influenced by the chaotic example $R^{\circ} \otimes R$ itself where we believe the antipode should be the switch isomorphism $(R^{\circ} \otimes R)^{\circ} \longrightarrow R^{\circ} \otimes R$. This is not a dualization in the sense of [DMS] but a dualization of the kind that arises in Barr's $*$-autonomous monoidal categories [Ba1].

Consequently we are led to study $*$-autonomy for enriched categories. In fact, we define $*$-autonomous promonoidal $\mathcal{V}$-categories and show this notion is preserved under convolution. There is always the canonical promonoidal structure on $\mathcal{A}^{op} \otimes \mathcal{A}$ (see the concluding remarks of [Da1]) which is $*$-autonomous (as remarked by Luigi Santocanale after the talk [Da4]) and leads under convolution to the tensor product of bimodules. The Chu construction as described in [Ba3] and [St4] is purely for ordinary categories: it needs the repetition and deletion of variables that are available in a cartesian closed base category such as Set. We vastly extend the notion of $*$-autonomy to include enriched categories and other contexts. We provide a general star-construction which leads to the Chu construction as a special case.

Equipped with this we can define when a Takeuchi $\times_R$-bialgebra is "Hopf". Then, by dualizing from k-algebras to k-coalgebras, we define quantum groupoids to be $*$-autonomous quantum categories.

## 2. Ordinary categories revisited

Let us consider Proposition 1.1 from a slightly different viewpoint. A left adjoint (or cocontinuous) functor $F : \mathrm{Set} / X \longrightarrow \mathrm{Set} / Y$ between slice categories is determined by its restriction to the elements $x : 1 \longrightarrow X$ of X, and so, by a functor

$$X \longrightarrow \mathrm{Set} / Y \xrightarrow{\sim} [Y, \mathrm{Set}],$$

where we regard the sets X and Y as discrete categories and write $[\mathcal{A}, \mathcal{B}]$ for the category of functors and natural transformations from $\mathcal{A}$ to $\mathcal{B}$. However, the functors $X \longrightarrow [Y, \mathrm{Set}]$ are in bijection with functors $S : X \times Y \longrightarrow \mathrm{Set}$ which we think of as matrices

$$S = \big( S(x;y) \big)_{(x,y) \in X \times Y} .$$

This gives us an equivalent (actually "biequivalent") way of looking at the 2-

category whose objects are (small) sets, whose morphisms $F : X \longrightarrow Y$ are cocontinuous functors $Set / X \longrightarrow Set / Y$, and whose 2-cells are natural transformations; however, rather than a 2-category we only have a bicategory which we call Mat(Set) (compare [BCSW] for example). Again, the objects are sets, the morphisms $S : X \longrightarrow Y$ are matrices, and the 2-cells $\theta : S \Rightarrow T$ are matrices of functions

$$\theta = \Big( \theta(x;y) : S(x;y) \longrightarrow T(x;y) \Big)_{(x,y) \in X \times Y} ;$$

vertical composition of 2-cells is defined by entrywise composition of functions, horizontal composition of morphisms $S : X \longrightarrow Y$ and $T : Y \longrightarrow Z$ is defined by matrix multiplication

$$(T \circ S)(x;z) = \sum_{y \in Y} S(x;y) \times T(y;z) ,$$

and horizontal composition is extended in the obvious way to 2-cells. We write $X : X \longrightarrow X$ for the identity matrix (or Kronecker delta):

$$X(x;y) = \begin{cases} 1 & \text{for } \ x = y, \\ \varnothing & \text{otherwise.} \end{cases}$$

Of course Mat(Set) is also biequivalent to the bicategory Span(Set) of spans (in the sense of Bénabou [Bé]) in the category Set of sets.

In fact, Mat(Set) is an autonomous monoidal bicategory in the sense of the authors [DS1]. That is, there is a reasonably well behaved tensor product pseudofunctor

$$\text{Mat(Set)} \times \text{Mat(Set)} \longrightarrow \text{Mat(Set)}$$

which is simply defined on objects by cartesian product of sets and likewise, by cartesian product entrywise, on morphisms and 2-cells. Each object $Y$ is actually self-dual since a matrix $X \times Y \longrightarrow Z$ can be identified with a matrix $X \longrightarrow Y \times Z$. This means that $Y \times Z$ is the internal hom in Mat(Set) of $Y$ and $Z$ (mimicking the fact that in finite-dimensional vector spaces the vector space of linear functions from $V$ to $W$ is isomorphic to $V^* \otimes W$). In particular, $X \times X$ is the internal endohom of $X$; and so we expect it to be a pseudomonoid in Mat(Set) (mimicking the fact that the internal endohom of an object in a monoidal category is an internal monoid).

Let us be more specific about this pseudomonoid structure on $X \times X$ in Mat(Set). The multiplication

$$P : (X \times X) \times (X \times X) \longrightarrow X \times X$$

is defined by $P(y_2, x_2, y_1, x_1 ; x, y) = X(y, x_1) \times X(y_1, x_2) \times X(y_2, x)$. The unit $J : 1 \longrightarrow X \times X$

is defined by $J(\bullet; x,y) = X(x;y)$. One easily checks the canonical associativity and unit isomorphisms

$$P \circ (P \times (X \times X)) \cong P \circ ((X \times X) \times P), \qquad P \circ (J \times (X \times X)) \cong X \times X \cong P \circ ((X \times X) \times J).$$

Thinking of the set $X \times X$ as a discrete category, we see that $P$, $J$ and these isomorphisms form a promonoidal structure on $X \times X$. Noting that, under the equivalence of categories

$$[X \times X, \mathrm{Set}] \xrightarrow{\sim} \mathrm{Set} / X \times X,$$

the convolution monoidal structure for $X \times X$ transports across the equivalence to the monoidal structure on $\mathrm{Set} / X \times X$ described in the Introduction, the following result becomes a corollary of Proposition 1.1.

**Proposition 2.1** *Categories with object set* $X$ *are equivalent to monoidal comonads on the internal endohom pseudomonoid* $X \times X$ *in the monoidal bicategory* Mat(Set).

It may be instructive to sketch a direct proof of this result. A monoidal comonad $G$ on $X \times X$ comes equipped with 2-cells

$$\delta : G \longrightarrow G \circ G, \;\; \varepsilon : G \longrightarrow X \times X, \;\; \mu : P \circ (G \times G) \longrightarrow G \circ P \;\text{ and }\; \eta : J \longrightarrow G \circ J,$$

subject to appropriate axioms. The mere existence of $\varepsilon$ is quite a strong condition since $X(x;u) \times X(y;v)$ is empty unless $x = u$ and $y = v$; so $G(x,y;u,v)$ is empty unless $x = u$ and $y = v$. This leads us to put

$$A(x,y) = G(x,y;x,y)$$

which defines the homsets of our category $A$. It is then easy to check that $\mu$ defines composition and $\eta$ provides the identities for the category $A$. We note finally that $\delta$ is forced to be a genuine diagonal morphism: we are dealing here with the categories of "commutative geometry".

## 3. Takeuchi bialgebroids

We are now ready to move from set theory to linear algebra. Let $k$ be any commutative ring and write $\mathcal{V}$ for the monoidal category of k-modules; we write $\otimes$ for the tensor product of k-modules. Monoids $R$ in $\mathcal{V}$ will be called k-algebras and we write $\mathcal{V}^R$ for the category of left R-modules. We recall the preliminaries of Morita theory. For k-algebras $R$ and $S$, a left adjoint (or cocontinuous) functor $F : \mathcal{V}^R \longrightarrow \mathcal{V}^S$ between module categories is, up to isomorphism, determined by its restriction to the full subcategory of $\mathcal{V}^R$ consisting of $R$ itself as a left R-module. This left S-module $F(R) = M$ is therefore actually a right R-, left S-bimodule which we call a *module from* $R$ *to* S and use the arrow notation $M : R \longrightarrow S$. (The fact that R is actually on the left of the arrow and S on the right, rather than the other way around, has to do with our convention to compose functions in the usual order.) We can also identify M with an object of $\mathcal{V}^{R^\circ \otimes S}$ where $R^\circ$ is the k-algebra R with opposite multiplication.

There is a 2-category whose objects are k-algebras, whose morphisms $R \longrightarrow S$ are left adjoint functors $F : \mathcal{V}^R \longrightarrow \mathcal{V}^S$, and whose 2-cells are natural transformations between such functors F; the compositions are the usual ones for functors and natural transformations. This 2-category is biequivalent to the bicategory $\mathrm{Mod}(\mathcal{V})$ whose objects are k-algebras, whose morphisms are modules $M : R \longrightarrow S$, and whose 2-cells are 2-sided module morphisms; the horizontal composite $N \circ M : R \longrightarrow T$ of $M : R \longrightarrow S$ and $N : S \longrightarrow T$ is the tensor product $N \otimes_S M$ of the modules M and N over S; vertical composition of 2-cells is the usual composition of module morphisms.

Indeed, like Mat(Set), the bicategory $\mathrm{Mod}(\mathcal{V})$ is autonomous monoidal. The tensor product is that of $\mathcal{V}$: k-algebras R and S are taken to the k-algebra $R \otimes S$, modules $M : R \longrightarrow S$ and $M' : R' \longrightarrow S'$ are taken to the module $M \otimes M' : R \otimes R' \longrightarrow S \otimes S'$, and module morphisms are tensored via the obvious inducement. The opposite k-algebra $S^\circ$ acts as a dual for S since the category of modules $R \otimes S \longrightarrow T$ is equivalent to the category of modules $R \longrightarrow S^\circ \otimes T$.

It follows that $R^\circ \otimes R$ is an internal endohom for R and, as such, is a pseudomonoid in $\mathrm{Mod}(\mathcal{V})$. The multiplication

$$P : (R^\circ \otimes R) \otimes (R^\circ \otimes R) \longrightarrow R^\circ \otimes R$$

is $P = R \otimes R \otimes R$ as a k-module, with the further actions defined by

$$(x \otimes y)(a \otimes b \otimes c)(x_1 \otimes y_1 \otimes x_2 \otimes y_2) = (yax_1) \otimes (y_1 b x_2) \otimes (y_2 c x)$$

for $a \otimes b \otimes c \in P$, $x \otimes y \in R^\circ \otimes R$ and $x_1 \otimes y_1 \otimes x_2 \otimes y_2 \in R^\circ \otimes R \otimes R^\circ \otimes R$. The unit

$$J : k \longrightarrow R^\circ \otimes R$$

is just $J = R$ as a k-module, with the furthur action $(x \otimes y)a = yax$. One easily checks that there are canonical isomorphisms

$$P \otimes_{R^e \otimes R^e} (R^e \otimes P) \cong P \otimes_{R^e \otimes R^e} (P \otimes R^e) \quad \text{and}$$

$$P \otimes_{R^e \otimes R^e} (R^e \otimes J) \cong R^e \cong P \otimes_{R^e \otimes R^e} (J \otimes R^e)$$

where we have used the traditional notation $R^e = R^\circ \otimes R$ for this pseudomonoid; the "e" superscript could be thought to stand for "endo" as well as the usual "envelope".

**Definition 3.1** A *Takeuchi bialgebroid* is a k-module R together with an opmonoidal monad on $R^e$ in the monoidal bicategory $\mathrm{Mod}(\mathcal{V})$.

To see that this definition agrees with that of $\times_R$-bialgebra as defined by Takeuchi [Tak] (and developed by [Lu], [Xu], [Sch], [BM] and [Szl]) we shall be more explicit about what an opmonoidal monad A on any pseudomonoid E involves.

In any monoidal bicategory $\mathcal{B}$ (with tensor product $\otimes$ and unit k ) we use the term *pseudomonoid* (or "monoidal object") for an object E equipped with a binary multiplication $P : E \otimes E \longrightarrow E$ and a unit $J : k \longrightarrow E$ which are associative and unital up to coherent invertible 2-cells. A *monoidal morphism* $f : E \longrightarrow E'$ is a morphism equipped with coherent 2-cells $P \circ (f \otimes f) \Rightarrow f \circ P$ and $J \Rightarrow f \circ J$. A *monoidal 2-cell* is one compatible with these last coherent 2-cells. With the obvious compositions, this defines a bicategory $\mathrm{Mon}\mathcal{B}$ of pseudomonoids in $\mathcal{B}$. For example, if $\mathcal{B}$ is the cartesian-monoidal 2-category Cat of categories, functors and natural transformations then $\mathrm{Mon}\mathcal{B}$ is the 2-category MonCat of monoidal categories, monoidal functors and monoidal natural transformations as in [EK].

We write $\mathcal{B}^{co}$ for the bicategory obtained from $\mathcal{B}$ on reversing 2-cells. We put

$$\mathrm{Opmon}\mathcal{B} = (\mathrm{Mon}\mathcal{B}^{co})^{co};$$

the objects are again pseudomonoids, the morphisms are *opmonoidal morphisms*, and the 2-cells are *opmonoidal 2-cells*. An *opmonoidal monad* in $\mathcal{B}$ is a monad in $\mathrm{Opmon}\mathcal{B}$.

A monoidal morphism $f : E \longrightarrow E'$ is called *strong* when the 2-cells $J \Rightarrow f \circ J$ and $P \circ (f \otimes f) \Rightarrow f \circ P$ are invertible. The inverses for these 2-cells equip such a strong f with the structure of opmonoidal morphism.

Now we return to the case of opmonoidal monads in $\mathcal{B} = \mathrm{Mod}(\mathcal{V})$. First of all, we have a module $A : E \longrightarrow E$. The monad structure consists of module morphisms

$$\mu : A \otimes_E A \longrightarrow A \qquad \text{and} \qquad \eta : E \longrightarrow A$$

satisfying the usual conditions of associativity and unitality:

$$\mu \circ (\mu \otimes_E 1_A) = \mu \circ (1_A \otimes_E \mu), \qquad \mu \circ (\eta \otimes_E 1_A) = 1_A = \mu \circ (1_A \otimes_E \eta).$$

The opmonoidal structure consists of module morphisms

$$\delta : A \otimes_E P \longrightarrow P \otimes_{E \otimes E} (A \otimes A) \quad \text{and} \quad \varepsilon : A \otimes_E J \longrightarrow J$$

satisfying the following conditions:

$$\begin{array}{ccc}
A \otimes_E P \otimes_{E^{\otimes 2}} (E \otimes P) \cong A \otimes_E P \otimes_{E^{\otimes 2}} (P \otimes E) & \xrightarrow{\delta \otimes 1} & P \otimes_{E^{\otimes 2}} A^{\otimes 2} \otimes_{E^{\otimes 2}} (P \otimes E) \\
 & & \wr\| \\
\Big\downarrow{\scriptstyle \delta \otimes 1} & & P \otimes_{E^{\otimes 2}} ((A \otimes_E P) \otimes A) \\
 & & \Big\downarrow{\scriptstyle 1 \otimes (\delta \otimes 1)} \\
P \otimes_{E^{\otimes 2}} A^{\otimes 2} \otimes_{E^{\otimes 2}} (E \otimes P) & & \\
\wr\| & & \\
P \otimes_{E^{\otimes 2}} (A \otimes (A \otimes_E P)) & \xrightarrow[1 \otimes (1 \otimes \delta)]{} & P \otimes_{E^{\otimes 2}} (A \otimes (P \otimes_{E^{\otimes 2}} A^{\otimes 2})) \cong P \otimes_{E^{\otimes 2}} ((P \otimes_{E^{\otimes 2}} A^{\otimes 2}) \otimes A)
\end{array}$$

$$\begin{array}{ccc}
A \otimes_E P \otimes_{E^{\otimes 2}} (E \otimes J) & \xrightarrow{\delta \otimes 1} & P \otimes_{E^{\otimes 2}} A^{\otimes 2} \otimes_{E^{\otimes 2}} (E \otimes J) \cong P \otimes_{E^{\otimes 2}} (A \otimes (A \otimes_E J)) \\
 & \searrow^{\cong} & \Big\downarrow{\scriptstyle 1 \otimes (1 \otimes \varepsilon)} \\
 & & A \cong P \otimes_{E^{\otimes 2}} (A \otimes J)
\end{array}$$

$$\begin{array}{ccc}
A \otimes_E P \otimes_{E^{\otimes 2}} (J \otimes E) & \xrightarrow{\delta \otimes 1} & P \otimes_{E^{\otimes 2}} A^{\otimes 2} \otimes_{E^{\otimes 2}} (J \otimes E) \cong P \otimes_{E^{\otimes 2}} ((A \otimes_E J) \otimes A) \\
 & \searrow^{\cong} & \Big\downarrow{\scriptstyle 1 \otimes (\varepsilon \otimes 1)} \\
 & & A \cong P \otimes_{E^{\otimes 2}} (J \otimes A)
\end{array}$$

$$\begin{array}{ccccc}
A \otimes_E A \otimes_E P & \xrightarrow{\mu \otimes 1} & A \otimes_E P & \xrightarrow{\delta} & P \otimes_{E^{\otimes 2}} A^{\otimes 2} \\
\Big\downarrow{\scriptstyle 1 \otimes \delta} & & & & \Big\uparrow{\scriptstyle 1 \otimes \mu^{\otimes 2}} \\
A \otimes_E P \otimes_{E^{\otimes 2}} A^{\otimes 2} & \xrightarrow{\delta \otimes 1} & P \otimes_{E^{\otimes 2}} A^{\otimes 2} \otimes_{E^{\otimes 2}} A^{\otimes 2} & \cong & P \otimes_{E^{\otimes 2}} (A \otimes_E A)^{\otimes 2}
\end{array}$$

$$\begin{array}{ccc}
A \otimes_E A \otimes_E J & \xrightarrow{\mu \otimes 1} & A \otimes_E J \\
\Big\downarrow{\scriptstyle 1 \otimes \varepsilon} & & \Big\downarrow{\scriptstyle \varepsilon} \\
A \otimes_E J & \xrightarrow[\varepsilon]{} & J
\end{array}
\qquad
\begin{array}{ccc}
P & \xrightarrow{\eta \otimes 1} & A \otimes_E P \\
 & \searrow^{1 \otimes \eta^{\otimes 2}} & \Big\downarrow{\scriptstyle \delta} \\
 & & P \otimes_{E^{\otimes 2}} A^{\otimes 2}
\end{array}
\qquad
\begin{array}{ccc}
J & \xrightarrow{\eta \otimes 1} & A \otimes_E J \\
 & \searrow^{1} & \Big\downarrow{\scriptstyle \varepsilon} \\
 & & J
\end{array}
\quad .$$

Notice in particular that $A$ becomes a k-algebra with multiplication defined by composing $\mu$ with the quotient morphism $A \otimes A \longrightarrow A \otimes_E A$ and with unit $\eta(1)$. Indeed, $\eta : E \longrightarrow A$ becomes a k-algebra morphism. Moreover, the structure on $A$ as a module $A : E \longrightarrow E$ is induced by $\eta : E \longrightarrow A$ via $e a e' = \eta(e)\, a\, \eta(e')$.

From time to time we will require special properties of bicategories such as Mod($\mathcal{V}$). In particular, at this moment, we need to point out that Mod($\mathcal{V}$) admits both the Kleisli and Eilenberg-Moore constructions for monads. For monads in 2-categories rather than bicategories, the universal nature of these constructions was defined in [St1]; however, for the kind of phenomenon for modules we are about to explain, a better reference is [St2].

First we recall that each k-algebra morphism $f : R \longrightarrow S$ leads to two modules $f_* : R \longrightarrow S$ and $f^* : S \longrightarrow R$ which are both equal to $S$ as k-modules but with the module actions defined by

$$s\,x\,r = s\,x\,f(r) \qquad \text{and} \qquad r\,y\,s = f(r)\,y\,s$$

for $x \in f_*$, $y \in f^*$, $r \in R$ and $s \in S$. What is more, there are module morphisms

$$R \longrightarrow f^* \otimes_S f_* \quad \text{and} \quad f_* \otimes_R f^* \longrightarrow S,$$

the former defined by $f$ and the latter defined by multiplication in $S$, forming the unit and counit of an adjunction in which $f^*$ is right adjoint to $f_*$.

Suppose $A : E \longrightarrow E$ is a monad on the k-algebra $E$ in the bicategory Mod($\mathcal{V}$). The multiplication $\mu : A \otimes_E A \longrightarrow A$ and unit $\eta : E \longrightarrow A$ morphisms compose with the quotient morphism $A \otimes A \longrightarrow A \otimes_E A$ and the unit $k \longrightarrow E$, respectively, provide the k-module $A$ with a k-algebra structure with $\eta : E \longrightarrow A$ becoming a morphism of k-algebras. Then $\mu$ can be regarded as a 2-cell

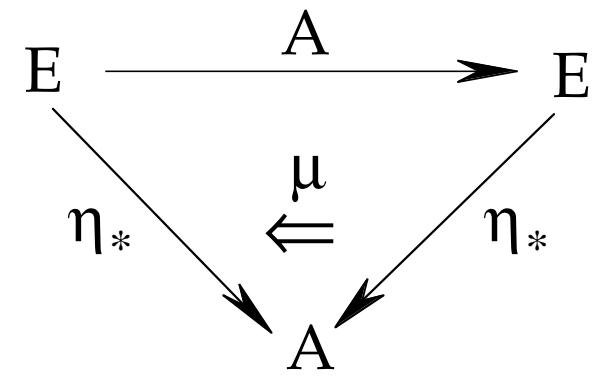

in Mod($\mathcal{V}$); it is a right action of the monoid $A$ on $\eta_*$. Indeed, this is the universal right action of $A$ on modules out of $E$; that is, the above triangle exhibits $A$ as the Kleisli construction for the monad $A$ on $E$. Since the homcategories of Mod($\mathcal{V}$) are cocomplete and composition with a given module preserves these colimits, the triangle

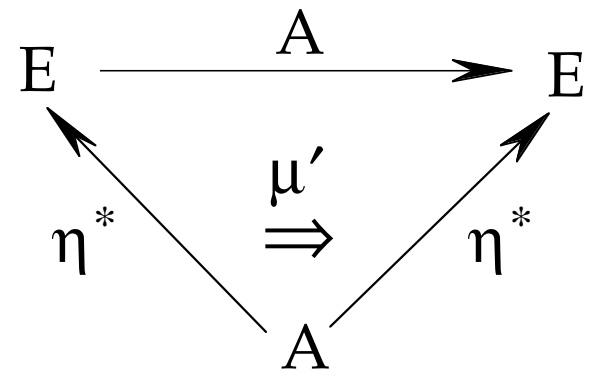

,

in which $\mu'$ is the mate of $\mu$ under the adjunction $\eta_* \dashv \eta^*$, exhibits $A$ as the Eilenberg-Moore construction for the monad A on E. That is, $\mu'$ is the universal left action of A of modules into E.

The following result abstracts Proposition 2.16 of [McC].

**Lemma 3.2** *If the monoidal bicategory $\mathcal{B}$ admits the Eilenberg-Moore construction for monads then so does* Opmon$\mathcal{B}$. *Furthermore, the forgetful morphism* Opmon$\mathcal{B} \longrightarrow \mathcal{B}$ *preserves the Eilenberg-Moore construction.*

In particular, this means that OpmonMod($\mathcal{V}$) admits the Eilenberg-Moore construction. (That the Kleisli construction exists for promonoidal monads was remarked in Section 3 of [Da2].)

**Proposition 3.3** *Suppose* E *is a pseudomonoid in* Mod($\mathcal{V}$) *and* $\eta : E \longrightarrow A$ *is a* k-*algebra morphism. There is an equivalence between the category of opmonoidal monad structures* $\mu$, $\delta$, $\varepsilon$ *on* $A : E \longrightarrow E$ *inducing* $\eta$ *and the category of pseudomonoid structures on* A *for which* $\eta^* : A \longrightarrow E$ *is a strong monoidal morphism.*

**Proof** In one direction, given the opmonoidal monad A on E inducing the given $\eta$, Lemma 3 lifts the triangle involving $\mu'$ to a triangle in OpmonMod($\mathcal{V}$) where it is again the Eilenberg-Moore construction. In particular, the adjunction $\eta_* \dashv \eta^*$ lifts to OpmonMod($\mathcal{V}$) and so, for general reasons explained in [Ke1], $\eta^* : A \longrightarrow E$ is strong monoidal. In the other direction, any k-algebra morphism $\eta : E \longrightarrow A$ always has the property that $\eta_*$ is opmonadic in Mod($\mathcal{V}$); that is, it supplies the Kleisli construction for the opmonoidal monad $\eta^* \otimes_A \eta_*$ on E generated by the adjunction $\eta_* \dashv \eta^*$. This opmonoidal monad has the form A, $\mu$, $\delta$, $\varepsilon$, $\eta$ as required. These two directions are the object functions for an obvious equivalence of categories. **QED**

It follows that a Takeuchi bialgebroid can equally be defined as consisting of a k-algebra R, a k-algebra morphism $\eta : R^e \longrightarrow A$, and a pseudomonoid structure on A for which $\eta^*$ is strong monoidal.

In preparation for interpreting Takeuchi bialgebroids in terms of module categories, we need to clarify further some monoidal terminology. The concepts are not new but the terminology is inconsistent in the literature.

We say that a monoidal $\mathcal{V}$-category $\mathcal{A}$ is *left closed* when, for all pairs of objects B,

C, there is an object $[B,C]_\ell$, called the *left internal hom of* B *and* C, for which there are isomorphisms

$$\mathcal{A}(A,[B,C]_\ell) \cong \mathcal{A}(A \otimes B, C),$$

$\mathcal{V}$-natural in A. A *right internal hom* $[B,C]_r$ satisfies

$$\mathcal{A}(A,[B,C]_r) \cong \mathcal{A}(B \otimes A, C).$$

We call a monoidal $\mathcal{V}$-category *closed* when it is both left and right closed. (This differs from Eilenberg-Kelly [EK] who use "closed" for left closed. However, they were mainly interested in the symmetric case where left closed implies right closed.)

As pointed out in [EK], if $\mathcal{A}$ and $\mathcal{X}$ are closed monoidal, a monoidal $\mathcal{V}$-functor $F : \mathcal{A} \longrightarrow \mathcal{X}$, with its (lax) constraints

$$\phi_0 : I \longrightarrow FI \quad \text{and} \quad \phi_{2;A,B} : FA \otimes FB \longrightarrow F(A \otimes B)$$

subject to axioms, could equally be called a *left closed* $\mathcal{V}$*-functor* since these constraints are in bijection with pairs

$$\phi_0 : I \longrightarrow FI \quad \text{and} \quad \phi^\ell_{2;B;C} : F[B,C]_\ell \longrightarrow [FB,FC]_\ell$$

satisfying corresponding axioms. Equally F could be called a *right closed* $\mathcal{V}$*-functor* since the constraints are in bijection with pairs

$$\phi_0 : I \longrightarrow FI \quad \text{and} \quad \phi^r_{2;A;C} : F[A,C]_r \longrightarrow [FA,FC]_r$$

satisfying corresponding axioms. We call a monoidal $\mathcal{V}$-functor F *normal* when $\phi_0$ is invertible. As usual we call F *strong monoidal* when it is normal and each $\phi_{2;A,B}$ is invertible. We define F to be *strong left closed* when it is normal and each $\phi^\ell_{2;B;C}$ is invertible; it is *strong right closed* when it is normal and each $\phi^r_{2;A;C}$ is invertible; and it is *strong closed* when it both strong left and strong right closed.

Pseudomonoid structures on A in $\mathrm{Mod}(\mathcal{V})$ are equivalent to closed monoidal structures on the $\mathcal{V}$-category $\mathcal{V}^A = \mathrm{Mod}(\mathcal{V})(k, A)$ of left A-modules; this is a special case of convolution in the sense of [Da1]. In fact, since k is a comonoid in $\mathrm{Mod}(\mathcal{V})$, we have a monoidal pseudofunctor

$$\mathrm{Mod}(\mathcal{V})(k,-) \; : \; \mathrm{Mod}(\mathcal{V}) \longrightarrow \mathcal{V} - \mathrm{Cat}\ ,$$

which, as such, takes pseudomonoids to pseudomonoids. Since it is representable by k, it also preserves Eilenberg-Moore constructions (and all weighted limits for that matter). This means that when we apply $\mathrm{Mod}(\mathcal{V})(k,-)$ to a Takeuchi bialgebroid $\eta : R^e \longrightarrow A$, we obtain a strong monoidal monadic functor

$$\mathcal{V}^{A} \longrightarrow \mathcal{V}^{R^e}.$$

Conversely, given a k-algebra morphism $\eta : R^e \longrightarrow A$, a $\mathcal{V}$-monoidal structure on $\mathcal{V}^A$, and a strong monoidal structure on the functor $\mathcal{V}^A \longrightarrow \mathcal{V}^{R^e}$, we obtain a Takeuchi bialgebroid structure on $\eta : R^e \longrightarrow A$. This is because $\mathcal{V}^A \longrightarrow \mathcal{V}^{R^e}$ has both adjoints and is conservative (= reflects isomorphisms), so is monadic; but being strong monoidal and colimit preserving, any monoidal structure on $\mathcal{V}^A$ will be automatically closed, reflecting the fact that the monoidal $\mathcal{V}$-category $\mathcal{V}^{R^e}$ is closed. Consequently, by [Da1], the monoidal structure on $\mathcal{V}^A$ is obtained by convolution of a pseudomonoid structure on A.

By the characterization Theorem 3.1 of [BM], we have shown that our Takeuchi bialgebroids are the $\times_R$-bialgebras of Takeuchi. We will see this in another way in the next section.

## 4. The lax monoidal operation $\times_R$

In order to define a bimonoid (or bialgebra) in a monoidal category, the monoidal category requires some kind of commutativity of the tensor product such as a braiding. A braiding can be regarded as a second monoidal structure on the category for which the new tensor is strongly monoidal with respect to the old. The so-called Eckmann-Hilton argument forces the new tensor to be isomorphic to the old and forces a braiding to appear (see [JS1]).

Ah, but what if the second tensor is only a lax multitensor and is only monoidal with respect to the old monoidal structure? Then there is certainly no need for the two structures to coincide. However, it is still possible to speak of a bimonoid: there is sufficient structure to express compatibility between a monoid structure for one tensor and a comonoid structure (on the same object) for the other tensor. After some preliminaries about right extensions in bicategories, we shall describe in detail just such a situation.

On top of the already discussed diverse properties and rich structure enjoyed by $\mathrm{Mod}(\mathcal{V})$, we also have the property that all right liftings and right extensions exist. Despite the terminology (from [St1] for example), these concepts are very familiar in the usual theory of modules.

Suppose M and M′ are modules $R \longrightarrow S$. We put

$$\mathrm{Hom}_R^S(M, M') = \mathrm{Mod}(\mathcal{V})(R,S)(M,M');$$

that is, traditionally, it is the k-module of left S, right R-bimodule morphisms from M

to $M'$. Now consider three modules as in the triangle

$$\begin{array}{ccc} R & \xrightarrow{\quad M \quad} & S \\ & {\scriptstyle N}\searrow \qquad \swarrow{\scriptstyle L} & \\ & T & \end{array}.$$

Let $\mathrm{Hom}_R(M,N) : S \longrightarrow T$ denote the k-module of right R-module morphisms with right S- and left T-actions defined by $(tfs)(m) = tf(sm)$ for

$$s \in S, \quad t \in T, \quad f \in \mathrm{Hom}_R(M,N) \quad \text{and} \quad m \in M.$$

Let $\mathrm{Hom}^T(L,N) : R \longrightarrow S$ denote the k-module of left T-module morphisms with right R- and left S-actions defined by $(sgr)(l) = g(ls)r$ for

$$r \in R, \quad s \in S, \quad g \in \mathrm{Hom}^T(L,N) \quad \text{and} \quad l \in L.$$

There are natural isomorphisms

$$\mathrm{Hom}_S^T(L,\mathrm{Hom}_R(M,N)) \cong \mathrm{Hom}_R^T(L \otimes_S M,N) \cong \mathrm{Hom}_R^S(M,\mathrm{Hom}^T(L,N)).$$

induced by evaluation morphisms

$$\mathrm{ev}_N^M : \mathrm{Hom}_R(M,N) \otimes_S M \longrightarrow N \quad \text{and} \quad \mathrm{ev}_N^L : L \otimes_S \mathrm{Hom}^T(L,N)) \longrightarrow N.$$

In bicategorical terms, $\mathrm{Hom}_R(M,N)$ is the right extension of $N$ along $M$, while $\mathrm{Hom}^T(L,N)$ is the right lifting of $N$ through $L$.

Consider any pseudomonoid $E$, with multiplication $P$ and unit $J$, in a monoidal bicategory $\mathcal{B}$ which admits all right extensions (where we have in mind $\mathcal{B} = \mathrm{Mod}(\mathcal{V})$). Then the endohom category $\mathrm{End}(E) = \mathcal{B}(E,E)$ becomes an oplax monoidal category ir the sense of [DS2] and [DS3], as follows. We define

$$P_n : E^{\otimes n} \longrightarrow E$$

to be the composite

$$E^{\otimes n} \xrightarrow{P \otimes E^{\otimes(n-2)}} E^{\otimes(n-1)} \xrightarrow{P \otimes E^{\otimes(n-3)}} \ldots \xrightarrow{P \otimes E} E^{\otimes 2} \xrightarrow{P} E$$

for $n \geq 2$, to be the identity of $n = 1$, and to be $J$ when $n = 0$. The coherence conditions for a pseudomonoid ensure that $P_m \cong P_n \circ (P_{m_1} \otimes \ldots \otimes P_{m_n})$ for each partition $\xi$ : $m_1 + \ldots + m_n = m$.

We define the multiple tensor $\underset{n}{\bullet}\big(M_1, \ldots, M_n\big)$ of objects $M_1, \ldots, M_n$ of $\mathrm{End}(E)$ to be the right extension of $P_n \circ \big(M_1 \otimes \ldots \otimes M_n\big)$ along $P_n$; that is,

$$\underset{n}{\bullet}\big(M_1, \ldots, M_n\big) = \mathrm{Hom}_{E^{\otimes n}}\Big(P_n, P_n \otimes_{E^{\otimes n}} \big(M_1 \otimes \ldots \otimes M_n\big)\Big).$$

The lax associativity constraint

$$\mu_\xi \ : \ \underset{n}{\bullet}\left(\underset{m_1}{\bullet}(M_{11},\ldots,M_{1m_1}),\ldots,\underset{m_n}{\bullet}(M_{n1},\ldots,M_{nm_n})\right)\longrightarrow \underset{m}{\bullet}(M_{11},\ldots,M_{nm_n})$$

for each partition $\xi : m_1+\ldots+m_n = m$ is, by using the right extension property of the target, induced by the morphism

$$\underset{n}{\bullet}\left(\underset{m_1}{\bullet}(M_{11},\ldots,M_{1m_1}),\ldots,\underset{m_n}{\bullet}(M_{n1},\ldots,M_{nm_n})\right)\circ P_m \longrightarrow P_m\circ(M_{11}\otimes\ldots\otimes M_{nm_n})$$

which, after "conjugation" with $P_m \cong P_n\circ(P_{m_1}\otimes\ldots\otimes P_{m_n})$, is the composite

$$\underset{n}{\bullet}\left(\underset{m_1}{\bullet}(M_{11},\ldots,M_{1m_1}),\ldots,\underset{m_n}{\bullet}(M_{n1},\ldots,M_{nm_n})\right)\circ P_n\circ(P_{m_1}\otimes\ldots\otimes P_{m_n})\xrightarrow{\ ev\circ 1\ }$$

$$P_n\circ\left(\underset{m_1}{\bullet}(M_{11},\ldots,M_{1m_1}),\ldots,\underset{m_n}{\bullet}(M_{n1},\ldots,M_{nm_n})\right)\circ(P_{m_1}\otimes\ldots\otimes P_{m_n})\xrightarrow{\ \cong\ }$$

$$P_n\circ\left(\left(\underset{m_1}{\bullet}(M_{11},\ldots,M_{1m_1})\circ P_{m_1}\right)\otimes\ldots\otimes\left(\underset{m_n}{\bullet}(M_{n1},\ldots,M_{nm_n})\circ P_{m_n}\right)\right)\xrightarrow{\ 1\circ(ev\otimes\ldots\otimes ev)\ }$$

$$P_n\circ\left(\left(P_{m_1}\circ\underset{m_1}{\bullet}(M_{11},\ldots,M_{1m_1})\right)\otimes\ldots\otimes\left(P_{m_n}\circ\underset{m_n}{\bullet}(M_{n1},\ldots,M_{nm_n})\right)\right)\xrightarrow{\ \cong\ }$$

$$P_n\circ(P_{m_1}\otimes\ldots\otimes P_{m_n})\circ(M_{11}\otimes\ldots\otimes M_{nm_n}).$$

The three axioms for a lax monoidal category can be verified. Since $P_1 : E\longrightarrow E$ is the identity, we see that $\underset{1}{\bullet} M = M$ ; so the lax monoidal structure on End(E) is normal.

As an endomorphism category End(E) is also a monoidal category for which the tensor product is composition. So End(E) is an object of the 2-category MonCat. Now MonCat is a monoidal 2-category with cartesian product as tensor. We will now see that End(E) is a lax monoid in MonCat.

**Proposition 4.1** *Regard* End(E) *as a monoidal category under composition. The functors* $\underset{n}{\bullet} : \mathrm{End}(E)^n\longrightarrow \mathrm{End}(E)$ *are equipped with canonical monoidal structures such that the substitutions* $\mu_\xi$ *are monoidal natural transformations.*

**Proof** The structure in question is the family of morphisms

$$\underset{n}{\bullet}(N_1,\ldots,N_n)\circ\underset{n}{\bullet}(M_1,\ldots,M_n) \longrightarrow \underset{n}{\bullet}(N_1\circ M_1,\ldots,N_n\circ M_n)$$

which, using the right extension property of the target, are induced by the composites

$$\underset{n}{\bullet}(N_1,\ldots,N_n)\circ\underset{n}{\bullet}(M_1,\ldots,M_n)\circ P_n\xrightarrow{\ 1\circ ev\ }$$

$$\underset{n}{\bullet}(N_1,\ldots,N_n)\circ P_n\circ(M_1\otimes\ldots\otimes M_n)\xrightarrow{\ ev\circ 1\ }$$

$$P_n \circ (N_1 \otimes \ldots \otimes N_n) \circ (M_1 \otimes \ldots \otimes M_n) \xrightarrow{\cong} P_n \circ \big((N_1 \circ M_1) \otimes \ldots \otimes (N_n \circ M_n)\big).$$

The compatibility of these morphisms with the lax associativity morphisms is readily verified. **QED**

A monoid for composition in End(E) is a monad on E in $\mathcal{B}$. We write MonEnd(E) for the category of monads on E; the morphisms are 2-cells between the endofunctors of the monads that are compatible with the units and multiplications. It follows from Proposition 5 that the lax monoidal structure on End(E) lifts to the category MonEnd(E).

The concept of comonoid makes sense in any lax monoidal category.

**Proposition 4.2** *A Takeuchi bialgebroid can equally be defined as a* k*-algebra* R *together with a comonoid in the lax monoidal category* MonEnd($R^e$).

**Proof** Both a Takeuchi bialgebroid $A : R^e \longrightarrow R^e$ and a comonoid in MonEnd($R^e$) start with a monad $A : R^e \longrightarrow R^e$ on $R^e$ in Mod($\mathcal{V}$). To make this a comonoid in MonEnd($R^e$) we need a comultiplication $\delta' : A \longrightarrow \underset{2}{\bullet}(A, A)$ and a counit $\varepsilon' : A \longrightarrow \underset{0}{\bullet}$ satisfying axioms. By the right extension properties of their targets, these morphisms are determined morphisms $\delta : A \circ P_2 \longrightarrow P_2 \circ (A \otimes A)$ and $\varepsilon : A \circ P_0 \longrightarrow P_0$, exactly as for a Takeuchi bialgebroid. The axioms which assert that $\delta'$ and $\varepsilon'$ should form a comonoid and should respect the monad structure translate precisely to the Takeuchi bialgebroid axioms (that is, opmonoidal monad axioms) on $\delta$ and $\varepsilon$. **QED**

The lax monoidal structure of Proposition 4.2 agrees with the operation $\times_R$ of Takeuchi [Tak].

## 5. Monoidal star autonomy

In this section we extend the theory of $*$–autonomous categories in the sense of Barr (see [Ba1], and, for the non-symmetric case, see [Ba3]) to enriched categories in the sense of Eilenberg-Kelly [EK]. The kind of duality present in a $*$–autonomous category is closer than compactness (also called rigidity or autonomy) to what is needed for an antipode in a bialgebroid or quantum category, and so for a concept of Hopf bialgebroid or quantum groupoid (see Example 7.4).

A $\mathcal{V}$-functor $F : \mathcal{A} \longrightarrow \mathcal{B}$ is called *eso* (for "essentially surjective on objects") when every object of $\mathcal{B}$ is isomorphic to one of the form FA for some object A of $\mathcal{A}$.

A *left star operation* for a monoidal $\mathcal{V}$-category $\mathcal{A}$ is an eso $\mathcal{V}$-functor

$$S_\ell :\ \mathcal{A} \longrightarrow \mathcal{A}^{\mathrm{op}}$$

together with a $\mathcal{V}$-natural family of isomorphisms (called the *left star constraint*)

$$\mathcal{A}(A \otimes B, S_\ell C) \cong \mathcal{A}(A, S_\ell(B \otimes C)).$$

It follows that $\mathcal{A}$ is left closed with $[B,C]_\ell \cong S_\ell(B \otimes D)$ where $S_\ell D \cong C$.

A *right star operation* for a monoidal $\mathcal{V}$-category $\mathcal{A}$ is an eso $\mathcal{V}$-functor

$$S_r :\ \mathcal{A}^{\mathrm{op}} \longrightarrow \mathcal{A}$$

together with a $\mathcal{V}$-natural family of isomorphisms (called the *right star constraint*)

$$\mathcal{A}(A \otimes B, S_r C) \cong \mathcal{A}(B, S_r(C \otimes A)).$$

It follows that $\mathcal{A}$ is then right closed with $[A,C]_r \cong S_r(E \otimes A)$ where $S_r E \cong C$.

A monoidal $\mathcal{V}$-category $\mathcal{A}$ is called $*$-*autonomous* when it is equipped with a left star operation which is fully faithful. Since it follows that $S_\ell$ is then an equivalence of $\mathcal{V}$-categories, we write $S_r$ for its adjoint equivalence so that the left star constraint can be written as

$$\mathcal{A}(A \otimes B, S_\ell C) \cong \mathcal{A}(B \otimes C, S_r A).$$

We see from this that $S_r$ is a right star operation and $*$-autonomy can equally be defined in terms of a fully faithful right star operation. It follows that $*$-autonomous monoidal $\mathcal{V}$-categories are closed, with internal homs given by the formulas

$$[B,C]_\ell \cong S_\ell(B \otimes S_r C) \quad \text{and} \quad [A,C]_r \cong S_r(S_\ell C \otimes A).$$

Notice that

$$\mathcal{A}(A, S_\ell I) \cong \mathcal{A}(I \otimes A, S_\ell I) \cong \mathcal{A}(A \otimes I, S_r I) \cong \mathcal{A}(A, S_r I),$$

so that $S_\ell I \cong S_r I$ (by the Yoneda Lemma). The object $S_\ell I$ is called the *dualizing object* and determines the left star operation via $[B, S_\ell I]_\ell \cong S_\ell B$.

For the reader interested in checking that our $*$-autonomous monoidal categories agree with Michael Barr's $*$-autonomous categories, we recommend Definition 2.3 of [Ba2] as the appropriate one for comparison. Also see [St3].

A monoidal category is autonomous if and only if there exists a left star operation $S_\ell$ and a family of $\mathcal{V}$-natural isomorphisms

$$S_\ell(A \otimes B) \cong S_\ell B \otimes S_\ell A.$$

If $\mathcal{A}$ is autonomous then taking the left dual provides a left star operation with isomorphism as required which *a fortiori* satisfy the conditions for a strong monoidal $\mathcal{V}$-functor. To see the less obvious implication, suppose we have an $S_\ell$ and the isomorphisms. Then $[B,C]_\ell \cong S_\ell(B \otimes D) \cong S_\ell D \otimes S_\ell B \cong C \otimes S_\ell B$ where $S_\ell D \cong C$, so $S_\ell B$ is a left dual for $B$. So every object $B$ has a left dual $S_\ell B$. However, every object $B$ is isomorphic to $S_\ell D$ for some D. This implies that D is a right dual for B.

## 6. Modules and promonoidal enriched categories

We shall discuss some basic facts about enriched categories and modules between them. Then we will review promonoidal categories and promonoidal functors in the enriched context. We obtain a result about restriction along a promonoidal functor.

Let $\mathcal{V}$ denotes any complete and cocomplete symmetric monoidal closed category. We write $\mathcal{V}$-Mod for the symmetric monoidal bicategory (in the sense of [DS1]) whose objects are $\mathcal{V}$-categories and whose hom-categories are defined by

$$\mathcal{V}\text{-Mod}(\mathcal{A}, \mathcal{B}) = \left[\mathcal{A}^{op} \otimes \mathcal{B}, \mathcal{V}\right].$$

The objects $M : \mathcal{A} \longrightarrow \mathcal{B}$ of $\mathcal{V}$-Mod($\mathcal{A}$ , $\mathcal{B}$) are called *modules from* $\mathcal{A}$ *to* $\mathcal{B}$. The composite of modules $M : \mathcal{A} \longrightarrow \mathcal{B}$ and $N : \mathcal{B} \longrightarrow \mathcal{C}$ is defined by the equation

$$(N \circ M)(A,C) = \int^{B} N(B,C) \otimes M(A,B);$$

the integral here is the "coend" in the sense of [DK] (also see [Ke2]). The tensor product for $\mathcal{V}$-Mod is the usual tensor product of $\mathcal{V}$-categories in the sense of [EK] (also see [Ke2]). Actually $\mathcal{V}$-Mod is autonomous since we have

$$\mathcal{V}\text{-Mod}\left(\mathcal{A} \otimes \mathcal{B}, \mathcal{C}\right) \cong \mathcal{V}\text{-Mod}\left(\mathcal{B}, \mathcal{A}^{op} \otimes \mathcal{C}\right)$$

since both sides are isomorphic to $\left[\mathcal{B}^{op} \otimes \mathcal{A}^{op} \otimes \mathcal{C}, \mathcal{V}\right]$.

We have reversed the direction of modules from that in [DS1–3] so that a promonoidal $\mathcal{V}$-category $\mathcal{A}$ is precisely a pseudomonoid (monoidal object) of $\mathcal{V}$-Mod (rather than $\mathcal{A}^{op}$ being such). The multiplication module $P : \mathcal{A} \otimes \mathcal{A} \longrightarrow \mathcal{A}$ and the unit module $J : \mathcal{I} \longrightarrow \mathcal{A}$ are equally $\mathcal{V}$-functors

$$P : \mathcal{A}^{op} \otimes \mathcal{A}^{op} \otimes \mathcal{A} \longrightarrow \mathcal{V} \quad \text{and} \quad J : \mathcal{A} \longrightarrow \mathcal{V},$$

and we have associativity constraints

$$\int^{X} P(X,C;D) \otimes P(A,B;X) \cong \int^{Y} P(A,Y;D) \otimes P(B,C;Y)$$

and unital constraints

$$\int^{X} P(X,A;B) \otimes JX \cong \mathcal{A}(A,B) \cong \int^{Y} P(A,Y;B) \otimes JY,$$

satisfying the usual two axioms (see [Da1]) which yield coherence. It is convenient to introduce the $\mathcal{V}$-functors

$$P_n : \underbrace{\mathcal{A}^{op} \otimes \ldots \otimes \mathcal{A}^{op} \otimes \mathcal{A}}_{n} \longrightarrow \mathcal{V},$$

for all natural numbers n, which we define as follows:

$$P_0 A = JA,\ \ P_1(A_1;A) = \mathcal{A}(A_1,A),\ \ P_2(A_1,A_2;A) = P(A_1,A_2;A) \text{ and}$$

$$P_{n+1}(A_1,\ldots,A_{n+1};A) = \int^X P(X,A_{n+1};A) \otimes P(A_1,\ldots,A_n;X).$$

We think of $P_n(A_1,\ldots,A_n;A)$ as the object of multimorphisms from $A_1,\ldots,A_n$ to A in $\mathcal{A}$. For example, when $\mathcal{A}$ is a monoidal $\mathcal{V}$-category, we have a promonoidal structure on $\mathcal{A}$ with

$$P_n(A_1,\ldots,A_n;A) \cong \mathcal{A}(A_1 \otimes \ldots \otimes A_n, A),$$

where the multitensor product is, say, bracketed from the left.

It will also be convenient to define a *multimorphism structure* on a $\mathcal{V}$-category $\mathcal{A}$ to be a sequence of $\mathcal{V}$-functors

$$P_n : \underbrace{\mathcal{A}^{op} \otimes \ldots \otimes \mathcal{A}^{op} \otimes \mathcal{A}}_{n} \longrightarrow \mathcal{V}$$

subject to no constraints. So a promonoidal structure is an example where all the $P_n$ are obtained from the particular ones for n = 0, 1, 2. A *multitensor structure* on $\mathcal{A}$ is a multimorphism structure for which each $P_n(A_1,\ldots,A_n;-)$ is representable; so we have objects $\bigotimes_n(A_1,\ldots,A_n)$ of $\mathcal{A}$ and a $\mathcal{V}$-natural family of isomorphisms

$$P_n(A_1,\ldots,A_n;A) \cong \mathcal{A}\Big(\bigotimes_n(A_1,\ldots,A_n), A\Big).$$

For example, when $\mathcal{A}$ is monoidal, we obtain $\bigotimes_n(A_1,\ldots,A_n)$ inductively from the cases n = 0, 1, and 2 where it is the unit, the identity functor, and the binary tensor product, respectively.

Suppose $\mathcal{A}$ and $\mathcal{X}$ are promonoidal $\mathcal{V}$-categories. A $\mathcal{V}$-functor $H : \mathcal{E} \longrightarrow \mathcal{A}$ is called *promonoidal* when it is equipped with $\mathcal{V}$-natural families of morphisms

$$\phi_{2;U,V;W} : P(U,V;W) \longrightarrow P(HU,HV;HW) \quad \text{and} \quad \phi_{0;U} : JU \longrightarrow JHU$$

that are compatible in the obvious way with the associativity and unital constraints. In fact, we can inductively define $\mathcal{V}$-natural families of morphisms

$$\phi_{n;U_1,\ldots,U_n;U} : P_n(U_1,\ldots,U_n;U) \longrightarrow P_n(HU_1,\ldots,HU_n;HU)$$

using the inductive definition of $P_n$. In particular, $\phi_{1;U;V} : \mathcal{E}(U,V) \longrightarrow \mathcal{A}(HU,HV)$ is the effect of H on homs. We say that H is *promonoidally fully faithful* when each $\phi_{n;U_1,\ldots,U_n;U}$ is invertible. We say F is *normal* when each $\phi_{0;U}$ is invertible.

A promonoidal $\mathcal{V}$-functor $H : \mathcal{E} \longrightarrow \mathcal{A}$ also gives rise in the obvious way to $\mathcal{V}$-natural families of morphisms

$$\bar{\phi}_{2;A,B;W} : \int^{U,V} P(U,V;W) \otimes \mathcal{A}(A,HU) \otimes \mathcal{A}(B,HV) \longrightarrow P(A,B;HW),$$

$$\phi^{\ell}_{2;U,B;C} : \int^{B,C} P(U,V;W) \otimes \mathcal{A}(B,HV) \otimes \mathcal{A}(HW,C) \longrightarrow P(HU,B;C), \text{ and}$$

$$\phi^{r}_{2;A,V;C} : \int^{U,W} P(U,V;W) \otimes \mathcal{A}(A,HU) \otimes \mathcal{A}(HW,C) \longrightarrow P(A,HV;C).$$

We need to say a little bit about convolution (see [Da1], [Da3] and [DS3]). For $\mathcal{V}$-categories $\mathcal{A}$ and $\mathcal{X}$ equipped with multimorphism structures, the *convolution multimorphism structure* on the $\mathcal{V}$-functor $\mathcal{V}$-category $[\mathcal{A},\mathcal{X}]$ is defined by

$$P_n(M_1,\dots,M_n;M) = \int_{A_1,\dots,A_n} [P_n(A_1,\dots,A_n;A), P_n(M_1A_1,\dots,M_nA_n;MA)]$$

whenever these ends all exist (for example, when $\mathcal{A}$ is small). In the case where $\mathcal{X}$ is multitensored, the convolution is also multitensored by the formula

$$\underset{n}{*}(M_1,\dots,M_n)(A) = \int^{A_1,\dots,A_n} P_n(A_1,\dots,A_n;A) \otimes \underset{n}{\bigotimes}(M_1A_1,\dots,M_nA_n),$$

provided the appropriate weighted colimits (expressed here by coends and tensors) exist in $\mathcal{X}$. In the case where $\mathcal{A}$ is promonoidal, if $\mathcal{X}$ is cocomplete closed monoidal then so is $[\mathcal{A},\mathcal{X}]$ (see [Da1]).

**Proposition 6.1** *Suppose* $H : \mathcal{E} \longrightarrow \mathcal{A}$ *is a normal promonoidal* $\mathcal{V}$*-functor. The restriction* $\mathcal{V}$*-functor*

$$[H,1] : [\mathcal{A},\mathcal{V}] \longrightarrow [\mathcal{E},\mathcal{V}]$$

*is a normal monoidal* $\mathcal{V}$*-functor. It is strong monoidal if and only if each* $\bar{\phi}_{2;A,B;W}$ *is invertible. It is strong left (respectively, strong right) closed if and only if each* $\phi^{\ell}_{2;U,B;C}$ *(respectively,* $\phi^{r}_{2;A,V;C}$*) is invertible.*

**Proof** The monoidal unital constraint for $[H,1]$ is $\phi_{0;U} : JU \longrightarrow JHU$. To obtain the associativity constraint, we use the Yoneda Lemma to replace

$$(MH * NH)W = \int^{U,V} P(U,V;W) \otimes MHU \otimes NHV$$

by the isomorphic expression

$$\int^{U,V,A,B} P(U,V;W) \otimes \mathcal{A}(A,HU) \otimes \mathcal{A}(B,HV) \otimes MHU \otimes NHV$$

and take the morphism into

$$(M * N)HW = \int^{A,B} P(A,B;HW) \otimes MA \otimes NB$$

of the form $\int^{A,B} \bar{\phi}_{2;A,B;W} \otimes 1 \otimes 1$ which is clearly invertible if $\bar{\phi}_{2;A,B;W}$ is. The converse comes by taking M and N to be representable and using Yoneda.

Similarly, the left closed constraint for $[H,1]$ is obtained by composing the morphism $\int_{B,C}\left[\phi^{\ell}_{2;U,B;C} \otimes 1, 1\right]$ from

$$[N,L]_{\ell} HU = \int_{B,C} [P(HU,B;C) \otimes NB, LC]$$

to

$$\int_{V,W,B,C} [P(U,V;W) \otimes \mathcal{A}(B,HV) \otimes \mathcal{A}(HW,C) \otimes NB, LC]$$

with the Yoneda isomorphism between this last expression and

$$[NH,LH]_{\ell} U = \int_{V,W} [P(U,V;W) \otimes NHV, LHW];$$

this constraint is clearly invertible if $\phi^{\ell}_{2;U,B;C}$ is, and the converse comes by taking N and L to be representable. The right closed case is dual. **QED**

## 7. Forms and promonoidal star autonomy

A problem with $*$–autonomy is that the common base categories (like the category of sets and the category of vector spaces) are not themselves $*$–autonomous. So we do not expect the convolution monoidal structure on $[\mathcal{A},\mathcal{V}]$ to be $*$–autonomous even when $\mathcal{A}$ is. We introduce the notion of *form* to address this problem: forms do exist on base categories and carry over to convolutions, while $*$–autonomy is to be equipped with a special kind of form. The definition of a $*$–autonomous promonoidal $\mathcal{V}$-category will be expressed in terms of forms.

A *form* for a promonoidal $\mathcal{V}$-category $\mathcal{A}$ is a module $\sigma : \mathcal{A} \otimes \mathcal{A} \longrightarrow I$ (where $I$ is the usual one-object $\mathcal{V}$-category) together with an isomorphism $\sigma \circ (P \otimes 1) \cong \sigma \circ (1 \otimes P)$. In other words, a form is a $\mathcal{V}$-functor

$$\sigma : \mathcal{A}^{op} \otimes \mathcal{A}^{op} \longrightarrow \mathcal{V}$$

together with a $\mathcal{V}$-natural family of isomorphisms

$$\int^{X} \sigma(X,C) \otimes P(A,B;X) \cong \int^{Y} \sigma(A,Y) \otimes P(B,C;Y)$$

called *form constraints.* Indeed, we can inductively obtain isomorphisms

$$\int^{X} \sigma(X, A_{n+1}) \otimes P_n(A_1,\dots,A_n;X) \cong \int^{Y} \sigma(A_1,Y) \otimes P_n(A_2,\dots,A_{n+1};Y)$$

called the *generalized form constraints.* A promonoidal $\mathcal{V}$-category with a chosen form is called *formal.*

For example, every object K of any promonoidal $\mathcal{V}$-category $\mathcal{A}$ defines a form $\sigma(A,B) = P(A,B;K)$; the form constraints are provided by the promonoidal associativity and unit constraints. Other examples are $*$-autonomous monoidal categories, as we shall soon discover. Moreover, we will also see that forms carry over to various constructions such as tensor products and general convolutions of $\mathcal{V}$-categories.

If $\mathcal{A}$ is monoidal, using Yoneda, the form constraints become

$$\sigma(A \otimes B, C) \cong \sigma(A, B \otimes C).$$

A form is called *continuous* when $\sigma(A,-)$ and $\sigma(-,B) : \mathcal{A}^{op} \longrightarrow \mathcal{V}$ are small (weighted) limit preserving for all objects A and B of $\mathcal{A}$.

**Proposition 7.1** *Let $\mathcal{A}$ and $\mathcal{X}$ be formal promonoidal $\mathcal{V}$-categories..*

(a) *If $\mathcal{A}$ and $\mathcal{X}$ are formal then the tensor product $\mathcal{A} \otimes \mathcal{X}$ with promonoidal structure*

$$P_n((A_1,X_1),\dots,(A_n,X_n);(A,X)) = P_n(A_1,\dots,A_n;A) \otimes P_n(X_1,\dots,X_n;X)$$

*admits the form* $\sigma((A,X),(B,Y)) = \sigma(A,B) \otimes \sigma(X,Y)$.

(b) *If $\mathcal{A}$ is small and $\mathcal{X}$ is cocomplete closed monoidal with a continuous form then the convolution monoidal $\mathcal{V}$-category $[\mathcal{A},\mathcal{X}]$ admits the continuous form*

$$\sigma(M,N) = \int_{A,B} [\sigma(A,B), \sigma(MA,NB)].$$

**Proof** (a) This is trivial.

(b) We have the calculation

$$\sigma(M * N, L) = \int_{U,C} [(M * N)U \otimes LC, \sigma(U,C)]$$

$$\cong \int_{U,C} \left[\sigma(U,C), \sigma\left(\int^{A,B} P(A,B;U) \otimes MA \otimes NB, LC\right)\right]$$

$$\cong \int_{U,A,B,C} [\sigma(U,C) \otimes P(A,B;U), \sigma(MA \otimes NB, LC)]$$

$$\cong \int_{U,A,B,C} [\sigma(A,U) \otimes P(B,C;U), \sigma(MA, NB \otimes LC)]$$

$$\cong \int_{U,A} \left[\sigma(A,U), \sigma\left(MA, \int^{B,C} P(B,C;U) \otimes NB \otimes LC\right)\right]$$

$$\cong \int_{U,A} [\sigma(A,U), \sigma(MA,(N*L)U)] \cong \sigma(M,N*L). \quad \textbf{QED}$$

A form $\sigma : \mathcal{A} \otimes \mathcal{A} \longrightarrow I$ transforms under the duality of $\mathcal{V}$-modules to a $\mathcal{V}$-module $\check{\sigma} : \mathcal{A} \longrightarrow \mathcal{A}^{op}$. We say the form $\sigma$ is *non-degenerate* when $\check{\sigma}$ is an equivalence as a $\mathcal{V}$-module (a Morita equivalence if you prefer). A form $\sigma$ is said to be *representable* when there exists a $\mathcal{V}$-functor $S_\ell : \mathcal{A} \longrightarrow \mathcal{A}^{op}$ and a $\mathcal{V}$-natural isomorphism

$$\sigma(A,B) \cong \mathcal{A}(A, S_\ell B).$$

A promonoidal $\mathcal{V}$-category is $*$-*autonomous* when it is equipped with a representable non-degenerate form. In fact, if $\mathcal{A}$ satisfies a minimal completeness condition ("Cauchy completeness") then "representable" is redundant. Notice that $S_\ell$ is necessarily an equivalence, with adjoint inverse $S_r$, say, and the form constraints have the cyclic appearance

$$P(A,B;S_\ell C) \cong P(B,C;S_r A).$$

More generally, using Yoneda, the generalized form constraints become

$$P_n(A_1,\ldots,A_n;S_\ell A_{n+1}) \cong \int^X \mathcal{A}(X,S_\ell A_{n+1}) \otimes P_n(A_1,\ldots,A_n;X)$$

$$\cong \int^X \sigma(X,A_{n+1}) \otimes P_n(A_1,\ldots,A_n;X) \cong \int^Y \sigma(A_1,Y) \otimes P_n(A_2,\ldots,A_{n+1};Y)$$

$$\cong \int^Y \mathcal{A}(Y,S_r A_1) \otimes P_n(A_2,\ldots,A_{n+1};Y) \cong P_n(A_2,\ldots,A_{n+1};S_r A_1).$$

A monoidal category is $*$-autonomous in the monoidal sense if and only if it is $*$-autonomous in the promonoidal sense.

**Corollary 7.2** *In Proposition 7.1, if* $\mathcal{A}$ *and* $\mathcal{X}$ *are* $*$*-autonomous then so are*

(a) $\mathcal{A} \otimes \mathcal{X}$ *and* (b) $[\mathcal{A}, \mathcal{X}]$.

**Proof** (a) $\sigma((A,X),(B,Y)) = \sigma(A,B) \otimes \sigma(X,Y) \cong \mathcal{A}(A,S_\ell B) \otimes \mathcal{X}(X,S_\ell Y)$
$\cong (\mathcal{A} \otimes \mathcal{X})((A,X),(S_\ell B,S_\ell Y)).$

(b) $\sigma(M,N) = \int_{A,B} [\sigma(A,B), \sigma(MA,NB)] \cong \int_{A,B} [\mathcal{A}(A,S_\ell B), \mathcal{X}(MA,S_\ell NB)]$
$\cong \int_B \mathcal{X}(MS_\ell B, S_\ell NB) \cong [\mathcal{A},\mathcal{X}](MS_\ell, S_\ell N) \cong [\mathcal{A},\mathcal{X}](M, S_\ell N S_r)$. **QED**

**Example 7.3** As noted in the final remarks of [Da1], for any $\mathcal{V}$-category $\mathcal{C}$, there is a canonical promonoidal structure on $\mathcal{C}^{op} \otimes \mathcal{C}$. It is explicitly defined by

$$P_0(C,D) = J(C,D) = \mathcal{C}(C,D) \quad \text{and}$$

$$P_2((D_1,C_1),(D_2,C_2);(C_3,D_3)) = \mathcal{C}(C_3,D_1)\otimes\mathcal{C}(C_1,D_2)\otimes\mathcal{C}(C_2,D_3).$$

More generally,

$$P_n((D_1,C_1),,\ldots,(D_n,C_n);(C_{n+1},D_{n+1})) = \mathcal{C}(C_{n+1},D_1)\otimes\mathcal{C}(C_1,D_2)\otimes\ldots\otimes\mathcal{C}(C_n,D_{n+1}).$$

After the lecture [Da4], Luigi Santocanale observed that $\mathcal{C}^{op}\otimes\mathcal{C}$ is $*$-autonomous. To be precise, define $S:(\mathcal{C}^{op}\otimes\mathcal{C})^{op}\longrightarrow\mathcal{C}^{op}\otimes\mathcal{C}$ by $S(D,C) = (C,D)$. Clearly

$$P_n((D_1,C_1),\ldots,(D_n,C_n);(C_{n+1},D_{n+1})) = P_n((D_2,C_2),\ldots,(D_{n+1},C_{n+1});(C_1,D_1)),$$

so that $S_r = S_\ell = S$ for $*$-autonomy. To relate this to our discussion of bialgebroids (Section 3), note that a k-algebra $\mathcal{C} = R$ is a one-object $\mathcal{V}$-category (for $\mathcal{V}$ the category of k-modules) and so the "chaotic bialgebroid" $\mathcal{C}^{op}\otimes\mathcal{C} = R^e$ is $*$-autonomous.

**Example 7.4** The notion of Hopf $\mathcal{V}$-algebroid as appearing as Definition 21 of [DS1] is an example of a $*$-autonomous promonoidal $\mathcal{V}$-category. Suppose that the $\mathcal{V}$-category $\mathcal{C}$ is *comonoidal* [Da1]; that is, $\mathcal{C}$ is a pseudomonoid (or monoidal object) in $(\mathcal{V}\text{-Cat})^{op}$: this means we have $\mathcal{V}$-functors $\Delta : \mathcal{C}\longrightarrow\mathcal{C}\otimes\mathcal{C}$ and $E : \mathcal{C}\longrightarrow\mathcal{I}$, coassociative and counital up to coherent $\mathcal{V}$-natural isomorphisms. It is easy to see that $\Delta$ must be given by the diagonal $\Delta C = (C,C)$ on objects. A multimorphism structure Q on $\mathcal{C}^{op}$ is then defined by

$$Q_n(C;C_1,\ldots,C_n) = \mathcal{C}(C,C_1)\otimes\ldots\otimes\mathcal{C}(C,C_n);$$

the actions on hom-objects require the $\mathcal{V}$-functors $\Delta$ and E. Indeed, Q defines a promonoidal structure (compare Section 5 of [Da1]). If this promonoidal $\mathcal{V}$-category is $*$-autonomous then the condition $Q(A,B;S_\ell C)\cong Q(B,C;S_r A)$ becomes

$$\mathcal{C}(A,S_\ell C)\otimes\mathcal{C}(B,S_\ell C)\cong\mathcal{C}(B,S_r A)\otimes\mathcal{C}(C,S_r A)\cong\mathcal{C}(B,S_r A)\otimes\mathcal{C}(A,S_\ell C),$$

which precisely gives the condition

$$\mathcal{C}(A,C)\otimes\mathcal{C}(B,C)\cong\mathcal{C}(B,S_r A)\otimes\mathcal{C}(A,C)$$

for the authors' concept of Hopf $\mathcal{V}$-algebroid.

A promonoidal functor $H : \mathcal{E}\longrightarrow\mathcal{A}$ between $*$-autonomous promonoidal $\mathcal{V}$-categories is called *$*$-autonomous* when it is equipped with a $\mathcal{V}$-natural transformation

$$\tau^\ell : HS_\ell\longrightarrow S_\ell H$$

such that the following diagram commutes

$$\begin{array}{ccccc}
P(U,V;S_\ell W) & \xrightarrow{\phi_{2;U,V;S_\ell W}} & P(HU,HV;HS_\ell W) & \xrightarrow{P(1,1;\tau^\ell)} & P(HU,HV;S_\ell HW) \\
\wr\| & & & & \wr\| \\
P(V,W;S_r U) & \xrightarrow[\phi_{2;V,W;S_r U}]{} & P(HV,HW;HS_r U) & \xrightarrow[P(1,1;\tau^r)]{} & P(HV,HW;S_r HU)
\end{array}$$

where $\tau^r : HS_r \longrightarrow S_r H$ is the mate of $\tau^\ell$ under the adjunction between $S_\ell$ and $S_r$. We call $H$ *strong* $*$-*autonomous* when $\tau^\ell$ is invertible; it follows that $\tau^r$ is invertible.

**Proposition 7.5** *Suppose* $H : \mathcal{E} \longrightarrow \mathcal{A}$ *is a strong* $*$-*autonomous promonoidal* $\mathcal{V}$-*functor. If the restriction* $\mathcal{V}$-*functor* $[H,1] : [\mathcal{A},\mathcal{V}] \longrightarrow [\mathcal{E},\mathcal{V}]$ *is strong monoidal then it is strong closed.*

**Proof** The idea of the proof is to use $*$-autonomy to cycle the criterion of Proposition 6.1 for $[H,1]$ to be strong monoidal into the criteria for it to be strong closed. The precise calculation for strong left closed is as follows:

$$\int^{V,W} P(U,V;W) \otimes \mathcal{A}(B,HV) \otimes \mathcal{A}(HW,C)$$

$$\cong \int^{V,W} P(U,V;W) \otimes \mathcal{A}(B,HV) \otimes \mathcal{A}(HW,S_\ell S_r C)$$

$$\cong \int^{V,W} P(U,V;W) \otimes \mathcal{A}(B,HV) \otimes \mathcal{A}(S_r C,S_r HW)$$

$$\cong \int^{V,W} P(U,V;W) \otimes \mathcal{A}(B,HV) \otimes \mathcal{A}(S_r C,HS_r W)$$

$$\cong \int^{V,W} P(U,V;S_\ell W) \otimes \mathcal{A}(B,HV) \otimes \mathcal{A}(S_r C,HW)$$

$$\cong \int^{V,W} P(V,W;S_r U) \otimes \mathcal{A}(B,HV) \otimes \mathcal{A}(S_r C,HW)$$

$$\cong P(B,S_r C;HS_r U) \cong P(B,S_r C;S_r HU) \cong P(HU,B;S_\ell S_r C) \cong P(HU,B;C).$$ **QED**

The next simple observation can be useful in this context.

**Proposition 7.6** *Suppose* $U : \mathcal{A} \longrightarrow \mathcal{X}$ *is any* $\mathcal{V}$-*functor with a left adjoint* $F$, *and suppose there are equivalences* $S : \mathcal{A} \longrightarrow \mathcal{A}^{op}$ *and* $S : \mathcal{X} \longrightarrow \mathcal{X}^{op}$ *such that* $S \circ U \cong U \circ S$. *Then* $U$ *has a right adjoint* $S^{-1} \circ F \circ S$ *and the monad* $T = U \circ F$ *generated by the original adjunction has a right adjoint comonad* $G = U \circ S^{-1} \circ F \circ S$. *Dually,* $F$ *has a left adjoint* $S^{-1} \circ U \circ S$. *A doubly infinite string of adjunctions is thereby created.*

**Proof** Clearly $U : \mathcal{A}^{op} \longrightarrow \mathcal{X}^{op}$ has $F$ as right adjoint whereas the mutually inverse equivalences $S$ and $S^{-1}$ are adjoint to each other on both sides. The results now follow by composing adjunctions. **QED**

An *opform* for a promonoidal $\mathcal{V}$-category $\mathcal{A}$ is a $\mathcal{V}$-functor

$$\sigma : \mathcal{A} \otimes \mathcal{A} \longrightarrow \mathcal{V}$$

and $\mathcal{V}$-natural isomorphisms

$$\int_X [P(A,B;X), \sigma(X,C)] \cong \int_Y [P(B,C;Y), \sigma(A,Y)],$$

called *opform constraints.* For a monoidal $\mathcal{V}$-category, we see by Yoneda's Lemma that an opform on $\mathcal{A}$ is the same as a form on $\mathcal{A}^{op}$. Moreover, in general, if $\sigma$ is a form on $\mathcal{A}$ and $K$ is any object of $\mathcal{V}$ then an opform $\sigma_K$ on $\mathcal{A}$ is defined by the equation

$$\sigma_K(A,B) = [\sigma(A,B), K].$$

**Proposition 7.7** *Let $\mathcal{A}$ be a small promonoidal $\mathcal{V}$-category. Each opform $\sigma$ for $\mathcal{A}$ determines a continuous form for the convolution monoidal* $\mathcal{V}$-category $[\mathcal{A}, \mathcal{V}]$ *via the formula*

$$\sigma(M,N) = \int_{A,B} [MA \otimes NB, \sigma(A,B)].$$

*Furthermore, every continuous form on $[\mathcal{A}, \mathcal{V}]$ arises thus from an opform on $\mathcal{A}$.*

**Proof** We have the calculation

$$\begin{aligned}
\sigma(M * N, L) &= \int_{U,C} [(M * N)U \otimes LC, \sigma(U,C)] \\
&\cong \int_{U,C} \left[ \int^{A,B} P(A,B;U) \otimes MA \otimes NB \otimes LC, \sigma(U,C) \right] \\
&\cong \int_{U,A,B,C} [MA \otimes NB \otimes LC, [P(A,B;U), \sigma(U,C)]] \\
&\cong \int_{U,A,B,C} [MA \otimes NB \otimes LC, [P(B,C;U), \sigma(A,U)]] \\
&\cong \int_{U,A,B,C} [MA, [P(B,C;U) \otimes NB \otimes LC, \sigma(A,U)]] \\
&\cong \int_{U,A} \left[ MA, \left[ \int^{B,C} P(B,C;U) \otimes NB \otimes LC, \sigma(A,U) \right] \right] \\
&\cong \sigma(M, N * L).
\end{aligned}$$

Conversely, any continuous form $\sigma$ on $[\mathcal{A}, \mathcal{V}]$ will have

$$\sigma(M,N) \cong \sigma\left(\int^{A} MA \otimes \mathcal{A}(A,-), \int^{B} NB \otimes \mathcal{A}(B,-)\right)$$

$$\cong \int_{A,B} \left[MA \otimes NB, \sigma(\mathcal{A}(A,-), \mathcal{A}(B,-))\right],$$

so that $\sigma$ will be determined by its value on representables. We define $\sigma$ for $\mathcal{A}$ by

$$\sigma(A,B) = \sigma(\mathcal{A}(A,-), \mathcal{A}(B,-)).$$

We have the calculation

$$\int_{U}\left[P(A,B;U), \sigma(\mathcal{A}(U,-), \mathcal{A}(C,-))\right] \cong \sigma\left(\int^{U} P(A,B;U) \otimes \mathcal{A}(U,-), \mathcal{A}(C,-)\right)$$

$$\cong \sigma(P(A,B;-), \mathcal{A}(C,-)) \cong \sigma(\mathcal{A}(A,-) * \mathcal{A}(B,-), \mathcal{A}(C,-)) \cong \sigma(\mathcal{A}(A,-), \mathcal{A}(B,-) * \mathcal{A}(C,-))$$

$$\cong \sigma(\mathcal{A}(A,-), P(B,C;-)) \cong \int_{V}\left[P(B,C;V), \sigma(\mathcal{A}(A,-), \mathcal{A}(V,-))\right].$$ **QED**

## 8. The star and Chu constructions

We adhere to the spirit of the review [St4] where the Chu construction is defined at the multimorphism level. The star construction on a multimorphism structure yields one that is $*$-autonomous. When applied to a promonoidal $\mathcal{V}$-category, the result may not be promonoidal — hence the need to work at the more general level.

For that, we define a general multimorphism structure to be *$*$-autonomous* when there exists an equivalence $S_\ell : \mathcal{A} \longrightarrow \mathcal{A}^{op}$ of $\mathcal{V}$-categories and a sequence of $\mathcal{V}$-natural isomorphisms

$$P_n(A_1,\dots,A_n;S_\ell A_{n+1}) \cong P_n(A_2,\dots,A_{n+1};S_r A_1)$$

where $S_r$ is an adjoint inverse for $S_\ell$.

In this section we will show how to modify a multimorphism structure, with a prescribed $S_\ell$, to obtain a $*$-autonomous one with the same $S_\ell$. We first need a natural definition: an *equivalence* $F : \mathcal{A} \longrightarrow \mathcal{B}$ of multimorphism structures is an equivalence $F$ of $\mathcal{V}$-categories together with natural isomorphisms

$$P_n(A_1,\dots,A_n;A) \cong P_n(FA_1,\dots,FA_n;FA);$$

the inverse equivalence of $F$ is obviously also a multimorphism equivalence.

Notice that, for any $*$-autonomous multimorphism structure, $S_\ell \circ S_\ell : \mathcal{A} \longrightarrow \mathcal{A}$ is a multimorphism equivalence: for we have the calculation

$$P_n(A_1,\dots,A_n;A) \cong P_n(A_1,\dots,A_n;S_r S_\ell A) \cong P_n(S_\ell A, A_1\dots,A_{n-1};S_\ell A_n)$$

$$\cong P_n(S_\ell A, A_1\dots,A_{n-1};S_r S_\ell S_\ell A_n) \cong P_n(S_\ell S_\ell A_n, S_\ell A, A_1,\dots,A_{n-2};S_\ell A_{n-1}) \cong$$

$$\dots \cong P_n(S_\ell S_\ell A_2,\dots,S_\ell S_\ell A_n, S_\ell A;S_\ell A_1) \cong P_n(S_\ell S_\ell A_1,\dots,S_\ell S_\ell A_n;S_\ell S_\ell A).$$

Now to our construction. Suppose we have a multimorphism structure $P$ on any $\mathcal{V}$-category $\mathcal{A}$ equipped with a contravariant $\mathcal{V}$-functor $S_\ell : \mathcal{A}^{op} \longrightarrow \mathcal{A}$ such that $S_\ell \circ S_\ell : \mathcal{A} \longrightarrow \mathcal{A}$ is an equivalence of multimorphism structures. It follows that $S_\ell$ is an equivalence; we write $S_r$ for the adjoint equivalence. The *starring* of this situation is the multimorphism structure $P^*$ on $\mathcal{A}$ defined by the formula

$$P^*_n(X_1, \ldots, X_n; S_\ell X_{n+1}) =$$

$$\int^{U_{ij}\,(1 \le i < j \le n+1)} \bigotimes_{m=1}^{n+1} P_n\big(U_{m\,m+1}, \ldots, U_{m\,n+1}, S_r U_{1\,m}, \ldots, S_r U_{m-1\,m}; S_r X_m\big).$$

**Proposition 8.1** *The starring* $P^*$ *produces a* $*$*-autonomous multimorphism structure on* $\mathcal{A}$ *with the given* $S_\ell$.

**Proof** Extend the definition of the $U_{ij}$ and $X_i$ by putting $U_{ji} = S_r U_{ij}$ and $X_{n+i+1} = S_r S_r X_i$. From the definition, we have

$$P^*_n(X_2, \ldots, X_{n+1}; S_r X_1) =$$

$$\int^{V_{ij}\,(1 \le i < j \le n+1)} \bigotimes_{m=1}^{n+1} P_n\big(V_{m\,m+1}, \ldots, V_{m\,n+1}, S_r V_{1\,m}, \ldots, S_r V_{m-1\,m}; S_r X_{m+1}\big),$$

which we notice is isomorphic to the formula for $P^*_n(X_1, \ldots, X_n; S_\ell X_{n+1})$ on making the change of variables $V_{ij} = U_{i+1\,j+1}$ and using the isomorphisms

$$P_n\big(U_{12}, \ldots, U_{1\,n+1}; S_r X_1\big) \cong P_n\big(S_r S_r U_{12}, \ldots, S_r S_r U_{1\,n+1}; S_r S_r S_r X_1\big).$$ **QED**

Let $C$ be a $\mathcal{V}$-category with a multimorphism structure $P$ and a multimorphism equivalence $T : C \longrightarrow C$. We suppose furthermore that $C$ is a comonoidal $\mathcal{V}$-category with derived promonoidal structure $Q$ as made explicit in Example 7.4. We require that $T : C^{op} \longrightarrow C^{op}$ is an equivalence for the multimorphism structure $Q$ (that is, that $T : C \longrightarrow C$ is a comonoidal equivalence).

We want to apply the star construction to $\mathcal{A} = C^{op} \otimes C$ with $S_\ell(C,D) = (D, T^{-1}C)$, so that $S_r(C,D) = (TD, C)$, and with the tensor product multimorphism structu $Q \otimes P$ for the $P$ and $Q$ as described in the last paragraph. Notice that $S_\ell S_\ell(C,D) = (T^{-1}C, T^{-1}C)$ so that $S_\ell \circ S_\ell : C^{op} \otimes C \longrightarrow C^{op} \otimes C$ is indeed a multimorphism equivalence.

Let us calculate the star $R^*$ of $R = Q \otimes P$:

$$R^*_n\big((X_1, Y_1), \ldots, (X_n, Y_n); (Y_{n+1}, T^{-1}X_{n+1})\big) =$$

$$\int^{(U_{ij},V_{ij})} \bigotimes_{m=1}^{n+1} \begin{pmatrix} Q_n\left(TY_m ; U_{m\,m+1}, \ldots, U_{m\,n+1}, TV_{1\,m}, \ldots, TV_{m-1\,m}\right) \\ \otimes P_n\left(V_{m\,m+1}, \ldots, V_{m\,n+1}, U_{1\,m}, \ldots, U_{m-1\,m}; X_m\right) \end{pmatrix}$$

$$\cong \int^{(U_{ij},V_{ij})} \bigotimes_{m=1}^{n+1} \begin{pmatrix} \mathcal{C}(TY_m, U_{m\,m+1}) \otimes \ldots \otimes \mathcal{C}(TY_m, U_{m\,n+1}) \otimes \mathcal{C}(Y_m, V_{1\,m}) \otimes \ldots \otimes \mathcal{C}(Y_m, V_{m-1\,m}) \\ \otimes P_n\left(V_{m\,m+1}, \ldots, V_{m\,n+1}, U_{1\,m}, \ldots, U_{m-1\,m}; X_m\right) \end{pmatrix}$$

$$\cong \int^{(U_{ij},V_{ij})} \bigotimes_{r<s}\left(\mathcal{C}(TY_r, U_{rs}) \otimes \mathcal{C}(Y_s, V_{rs})\right) \otimes \bigotimes_{m=1}^{n+1} P_n\left(V_{m\,m+1}, \ldots, V_{m\,n+1}, U_{1\,m}, \ldots, U_{m-1\,m}; X_m\right)$$

$$\cong \bigotimes_{m=1}^{n+1} P_n\left(Y_{m+1}, \ldots, Y_{n+1}, TY_1, \ldots, TY_{m-1}; X_m\right)$$

which has the same shape as the multimorphism structure described in [St4].

**Proposition 8.2** *In the situation just described, if* P *is actually a monoidal structure on* $\mathcal{C}$*, then* $R^*$ *is a* $*$*-autonomous promonoidal structure on* $\mathcal{C}^{op} \otimes \mathcal{C}$.

**Proof** After Proposition 8.1, it suffices to show that $R^*$ is promonoidal. We need to see that each $R^*_n$ is determined by the $n = 0$ and $n = 2$ cases. The general calculation is by induction so we trust that the following exemplary step will be sufficient indication for the reader:

$$\int^{A_1,B_1} R^*_2\left((X_1,Y_1),(X_2,Y_2);(B_1,T^{-1}A_1)\right) \otimes R^*_2\left((TB_1,A_1),(X_3,Y_3);(B_4,T^{-1}A_4)\right)$$

$$\cong \int^{A_1,B_1} \begin{pmatrix} P_2(Y_2,B_1;X_1) \otimes P_2(B_1,TY_1;X_2) \otimes P_2(TY_1,TY_2;A_1) \\ \otimes P_2(Y_3,Y_4;TB_1) \otimes P_2(Y_4,A_1;X_3) \otimes P_2(A_1,TY_3;X_4) \end{pmatrix}$$

$$\cong \int^{A_1,B_1} \begin{pmatrix} \mathcal{C}(Y_2 \otimes B_1,X_1) \otimes \mathcal{C}(B_1 \otimes TY_1,X_2) \otimes \mathcal{C}(TY_1 \otimes TY_2,A_1) \\ \otimes \mathcal{C}(Y_3 \otimes Y_4,TB_1) \otimes \mathcal{C}(Y_4 \otimes A_1,X_3) \otimes \mathcal{C}(A_1 \otimes TY_3,X_4) \end{pmatrix}$$

$$\cong \begin{pmatrix} \mathcal{C}(Y_2 \otimes T(Y_3 \otimes Y_4),X_1) \otimes \mathcal{C}(T(Y_3 \otimes Y_4) \otimes TY_1,X_2) \\ \otimes \mathcal{C}(Y_4 \otimes TY_1 \otimes TY_2,X_3) \otimes \mathcal{C}(TY_1 \otimes TY_2 \otimes TY_3,X_4) \end{pmatrix}$$

$$\cong P_3(Y_2,TY_3,TY_4;X_1) \otimes P_3(TY_3,TY_4,TY_1;X_2) \otimes P_3(Y_4,TY_1,TY_2;X_3) \otimes P_3(TY_1,TY_2,TY_3;X_4)$$

$$\cong R^*_3\left((X_1,Y_1),(X_2,Y_2),(X_3,Y_3);(Y_4,T^{-1}X_4)\right).$$ **QED**

**Proposition 8.3** *In the situation of the Proposition 8.2, further suppose that* P *is closed monoidal and that the comonoidal structure on* $\mathcal{C}$ *is representable by an object* K*, an operation* $B \bullet C$*, and* $\mathcal{V}$*-natural isomorphisms*

$$\mathcal{C}(A,K) \cong I \quad \text{and} \quad \mathcal{C}(A,B \bullet C) \cong \mathcal{C}(A,B) \otimes \mathcal{C}(A,C)$$

*where the right-hand sides require the counit and comultiplication for their effects on homs. Then* $R^*$ *is a* $*$*-autonomous monoidal structure on* $\mathcal{C}^{op} \otimes \mathcal{C}$.

**Proof** We have the calculations:

$$\begin{aligned}
& R_2^*\big((X_1, Y_1), (X_2, Y_2); (Y_3, T^{-1}X_3)\big) \cong P_2(Y_2, Y_3; X_1) \otimes P_2(Y_3, TY_1; X_2) \otimes P_2(TY_1, TY_2; X_3) \\
& \cong \mathcal{C}(Y_2 \otimes Y_3, X_1) \otimes \mathcal{C}(Y_3 \otimes TY_1, X_2) \otimes \mathcal{C}(TY_1 \otimes TY_2, X_3) \\
& \cong \mathcal{C}\big(Y_3, [Y_2, X_1]_r\big) \otimes \mathcal{C}\big(Y_3, [TY_1, X_2]_\ell\big) \otimes \mathcal{C}(TY_1 \otimes TY_2, X_3) \\
& \cong \mathcal{C}\big(Y_3, [Y_2, X_1]_r \bullet [TY_1, X_2]_\ell\big) \otimes \mathcal{C}(Y_1 \otimes Y_2, T^{-1}X_3) \\
& \cong \big(\mathcal{C}^{op} \otimes \mathcal{C}\big)\Big(\big([Y_2, X_1]_r \bullet [TY_1, X_2]_\ell, Y_1 \otimes Y_2\big), (Y_3, T^{-1}X_3)\Big)
\end{aligned}$$

and

$$R_0^*(Y, T^{-1}X) \cong P_0(X) \cong \mathcal{C}(I, X) \cong \mathcal{C}(Y, K) \otimes \mathcal{C}(T^{-1}I, T^{-1}X) \cong \big(\mathcal{C}^{op} \otimes \mathcal{C}\big)\big((K, I), (Y, T^{-1}X)\big),$$

so that $\mathcal{C}^{op} \otimes \mathcal{C}$ is monoidal with unit $(K, I)$ and tensor product

$$(X_1, Y_1) \otimes (X_2, Y_2) = \big([Y_2, X_1]_r \bullet [TY_1, X_2]_\ell, Y_1 \otimes Y_2\big).$$ **QED**

A particular case of Proposition 8.3 is the Chu construction of [Ba3]. Here $\mathcal{V}$ is the category of sets with cartesian monoidal structure (although any cartesian closed base would do). Then every $\mathcal{V}$-category $\mathcal{C}$ is comonoidal via the diagonal functor $\Delta$. The representability of this structure as required in Proposition 8.3 amounts to $\mathcal{C}$ having finite limits; so $K$ is the terminal object and $B \bullet C = B \times C$ is the product of $B$ and $C$. Then $R^*$ is the $*$-autonomous monoidal structure on $\mathcal{C}^{op} \otimes \mathcal{C}$ arising from any monoidal closed category $\mathcal{C}$ with finite products and a monoidal endoequivalence $T$.

However, the case of finite products for ordinary categories is not the only example where the representable comonoidal structure can be found. For any $\mathcal{V}$, such structure exists for example on any $\mathcal{C}$ which is a free $\mathcal{V}$-category on an ordinary category with finite products.

## 9. Star autonomy in monoidal bicategories

In order to exploit duality, we need to generalise the notion of star autonomy to pseudomonoids in a monoidal bicategory $\mathcal{B}$. The work of Sections 5 to 8 is a special case taking place in the autonomous monoidal bicategory $\mathcal{V}$-Mod of $\mathcal{V}$-categories and $\mathcal{V}$-modules as defined in Section 6.

As mentioned in Section 3, for pseudomonoids $A$ and $E$ in $\mathcal{B}$, where we write p and j for the multiplications and units, a *monoidal morphism* $g : A \longrightarrow E$ is a morphism equipped with coherent 2-cells

$$\phi_2 : p \circ (g \otimes g) \Rightarrow g \circ p \text{ and } \phi_0 : j \Rightarrow g \circ j.$$

The morphism is called *strong monoidal* when $\phi_2$ and $\phi_0$ are invertible. When $g$ has a left adjoint $h$, there are 2-cells

$$\phi_2^\ell : h \circ p \circ (1 \otimes g) \Rightarrow p \circ (h \otimes 1) \quad \text{and} \quad \phi_2^r : h \circ p \circ (g \otimes 1) \Rightarrow p \circ (1 \otimes h)$$

obtained from $\phi_2$ as mates under adjunction. We say $g$ is *strong left [right] closed* when $\phi_2^\ell$ [respectively, $\phi_2^r$ ] is invertible; it is strong closed when it is both.

For a pseudomonoid $A$ in $\mathcal{B}$, the category $\mathcal{B}(I, A)$ is monoidal with tensor product defined by

$$m * n = p \circ (m \otimes n).$$

The internal homs, provided $\mathcal{B}$ has the relevant right liftings, are defined as follows: $[n, r]_\ell$ is the right lifting of $r$ through $p \circ (1_A \otimes n)$ while $[m, r]_r$ is the right lifting of $r$ through $p \circ (m \otimes 1_A)$.

**Proposition 9.1** *If* $g : A \longrightarrow E$ *is a strong monoidal morphism between pseudomonoids then* $\mathcal{B}(I, g) : \mathcal{B}(I, A) \longrightarrow \mathcal{B}(I, E)$ *is a strong monoidal functor. If* $g$ *has a left adjoint* $h$ *and is strong closed then the functor* $\mathcal{B}(I, g)$ *is strong closed.*

**Proof** For the first sentence we have

$$\begin{aligned}\mathcal{B}(I, g)(m * n) &= g \circ p \circ (m \otimes n) \cong p \circ (g \otimes g) \circ (m \otimes n) \\ &\cong p \circ (g \circ m) \otimes (g \circ n) \cong (g \circ m) * (g \circ n) \\ &\cong \mathcal{B}(I, g)(m) * \mathcal{B}(I, g)(n).\end{aligned}$$

For the second sentence consider the diagram

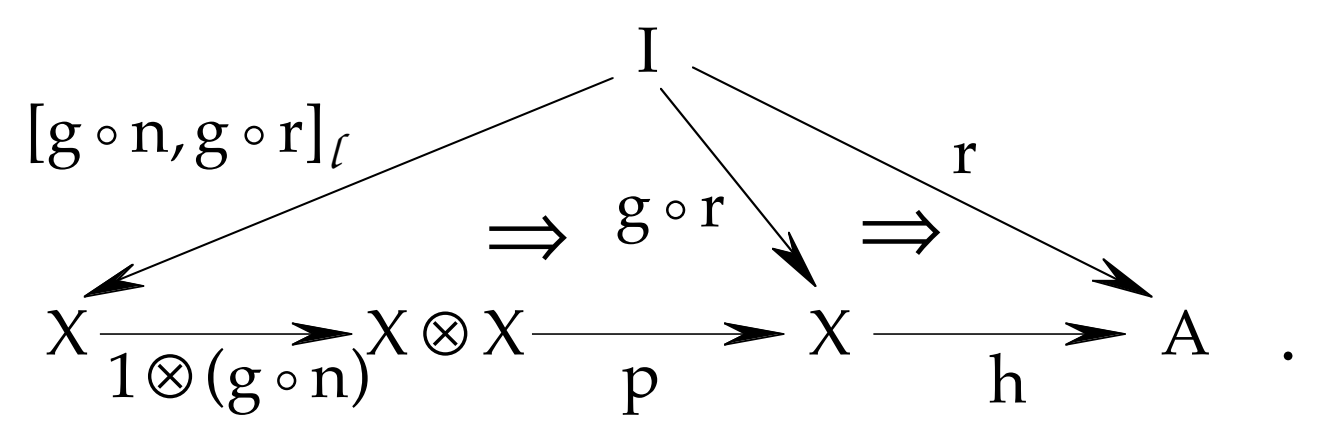

The right-hand triangle is a right lifting since $h$ is left adjoint to $g$. The left-hand triangle is a right lifting by definition of the left internal hom. So the outside triangle exhibits $[g \circ n, g \circ r]_\ell$ as a right lifting of $r$ along the bottom composite. However, if $g$ is strong left closed, the bottom composite is isomorphic to

$$h \circ p \circ (1 \otimes g) \circ (1 \otimes n) \cong p \circ (h \otimes 1) \circ (1 \otimes n) \cong p \circ (1 \otimes n) \circ h.$$

However, the right lifting of $r$ through $p \circ (1 \otimes n)$ is $[n, r]_\ell$, and the right lifting of $[n,r]_\ell$ through $h$ is $g \circ [n,r]_\ell$. So we have $g \circ [n,r]_\ell \cong [g \circ n, g \circ r]_\ell$ proving $\mathcal{B}(I,g)$ strong left closed. Right closedness is dual. **QED**

A *form* for a pseudomonoid $A$ in $\mathcal{B}$ is a morphism $\sigma : A \otimes A \longrightarrow I$ together with an isomorphism

$$\begin{array}{ccc} A^{\otimes 3} & \xrightarrow{\;p \otimes 1_A\;} & A^{\otimes 2} \\ {\scriptstyle 1_A \otimes p}\downarrow & \overset{\gamma}{\cong} & \downarrow{\scriptstyle \sigma} \\ A^{\otimes 2} & \xrightarrow[\;\sigma\;]{} & I \end{array}$$

called the *form constraint.*

In the bicategories $\mathcal{B}$ that we have in mind there are special morphisms (as abstracted by Wood [Wo]). The special morphisms $h$ have right adjoints $h^*$ and, in some cases, are precisely the morphisms with right adjoints, sometimes called *maps* in $\mathcal{B}$. For example, in Mat(Set) the maps are precisely the matrices arising from functions, and these are the special morphisms we want. For $\mathrm{Mod}(\mathcal{V})$, the special morphisms are those modules isomorphic to $h_*$ for some algebra morphism $h$. In the bicategory of $\mathcal{V}$-categories and $\mathcal{V}$-modules the special modules are those arising from $\mathcal{V}$-functors.

Suppose $\mathcal{B}$ has selected special maps and that $\mathcal{B}$ is autonomous. Each form $\sigma : A \otimes A \longrightarrow I$ corresponds to a morphism $\hat{\sigma} : A \longrightarrow A^{\circ}$. We say that the form $\sigma$ is *representable* when $\hat{\sigma}$ is isomorphic to a special map. We say that $\sigma$ is *non-degenerate* when $\hat{\sigma}$ is an equivalence.

A pseudomonoid in $\mathcal{B}$ is defined to be $*$*-autonomous* when it is equipped with a non-degenerate representable form.

An opmorphism $h : E \longrightarrow A$ between $*$-autonomous pseudomonoids is called $*$-*autonomous* when there is an isomorphism

$$\begin{array}{ccc} E \otimes E & \xrightarrow{\quad\sigma\quad} & I \\ {\scriptstyle h \otimes h}\downarrow & \overset{\tau}{\Leftarrow} & \\ A \otimes A & \xrightarrow[\quad\sigma\quad]{} & I \end{array}$$

such that the following equation holds:

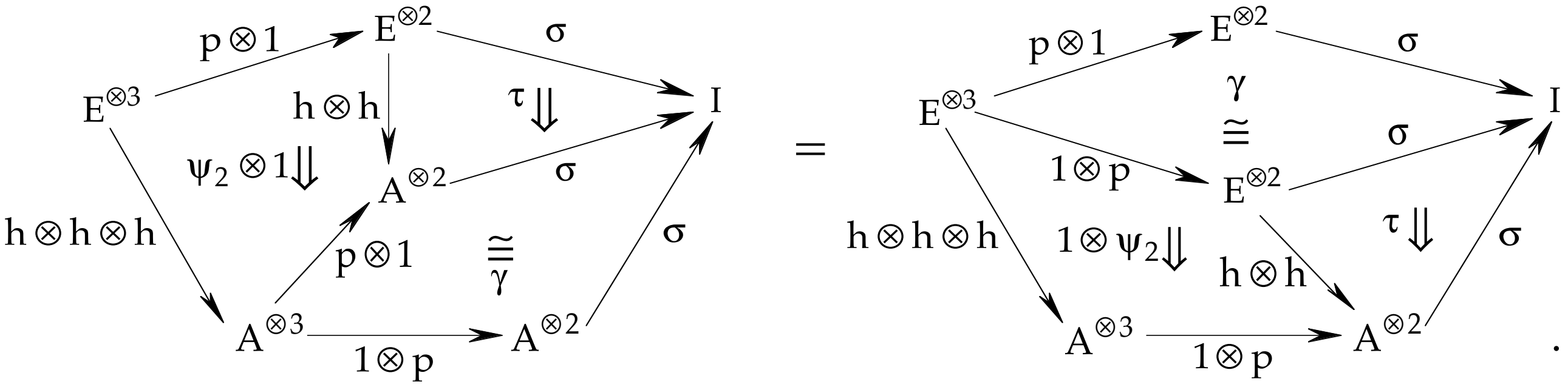

We are particularly interested in opmorphisms $h$ that are maps. Then the right adjoint $h^*$ is a morphism of pseudomonoids. Under these circumstances we define $h$ to be *strong $*$-autonomous* when the mate

$$\tau^{\ell} \; : \; \sigma \circ (h^* \otimes 1) \Rightarrow \sigma \circ (1 \otimes h)$$

of $\tau$ is invertible. It follows that $\tau^r \; : \; \sigma \circ (1 \otimes h^*) \Rightarrow \sigma \circ (h \otimes 1)$ is also invertible.

**Proposition 9.2** *Suppose* $h : E \longrightarrow A$ *is a strong $*$-autonomous special opmorphism between $*$-autonomous pseudomonoids in $\mathcal{B}$. If* $h^*$ *is strong monoidal then* $h^*$ *is strong closed.*

**Proof** We have the calculation

$$\sigma \circ (h \otimes 1) \circ (p \otimes 1) \circ (1 \otimes h^* \otimes 1) \;\cong\; \sigma \circ (1 \otimes h^*) \circ (p \otimes 1) \circ (1 \otimes h^* \otimes 1)$$

$$\cong\; \sigma \circ (p \otimes 1) \circ (1 \otimes h^*) \circ (1 \otimes h^* \otimes 1) \;\cong\; \sigma \circ (1 \otimes p) \circ (1 \otimes h^* \otimes h^*)$$

$$\cong\; \sigma \circ (1 \otimes h^*) \circ (1 \otimes p) \;\cong\; \sigma \circ (h \otimes 1) \circ (1 \otimes p) \;\cong\; \sigma \circ (1 \otimes p) \circ (h \otimes 1 \otimes 1)$$

$$\cong\; \sigma \circ (1 \otimes p) \circ (h \otimes 1 \otimes 1) .$$

It follows that $\hat{\sigma} \circ h \circ p \circ (1 \otimes h^*) \;\cong\; \hat{\sigma} \circ p \circ (h \otimes 1)$. Left strong closedness follows since $\sigma$ is non-degenerate. Right closedness is dual. **QED**

Motivated by Proposition 3.3, we define *basic data* in an autonomous monoidal bicategory $\mathcal{B}$ to consist of an object $R$ equipped with a special opmorphism $h : R^{\circ} \otimes R \longrightarrow A$ into a pseudomonoid $A$ such that $h^*$ is strong monoidal. Here $R^{\circ} \otimes R$ has the canonical endohom pseudomonoid structure. The basic data is called *Hopf* when $A$ is equipped with a $*$-autonomous structure and $h$ is strong $*$-autonomous.

From basic data, by applying the pseudofunctor $\mathcal{B}(I,-) : \mathcal{B} \longrightarrow \mathrm{Cat}$, we obtain an adjunction

$$\mathcal{B}(I,h) \dashv \mathcal{B}(I,h^*) \; : \; \mathcal{B}(I,A) \longrightarrow \mathcal{B}(I,R^{\circ} \otimes R)$$

which transports via the equivalence $\mathcal{B}(I,R^{\circ} \otimes R) \xrightarrow{\sim} \mathcal{B}(R,R)$ to an adjunction between $\mathcal{B}(I,A)$ and $\mathcal{B}(R,R)$. The pseudomonoid structure on A induces a monoidal structure on $\mathcal{B}(I,A)$ and the canonical endohom pseudomonoidal structure on $R^\circ \otimes R$ induces the monoidal structure on $\mathcal{B}(R,R)$ whose tensor product is composition in $\mathcal{B}$. Since $h^*$ is strong monoidal, the right adjoint $\mathcal{B}(I,A) \longrightarrow \mathcal{B}(R,R)$ is strong monoidal. By Propositions 9.1 and 9.2, this right adjoint is also strong closed in the Hopf case.

## 10. Ordinary groupoids revisited

Let us return to the definition of ordinary category as formulated in Propositions 1.1 and 2.1. Let G be a monoidal comonad on the internal endohom pseudomonoid $X \times X$ in the monoidal bicategory Mat(Set). Recall that $G(x,y;u,v)$ is empty unless $x = u$ and $y = v$, and we put

$$\mathbf{A}(x,y) = G(x,y;x,y)$$

which defines the homsets of our category $\mathbf{A}$. Let A denote the set of arrows of the category $\mathbf{A}$ ; we have the triangle which is the universal coaction of G on a morphism into $X \times X$; it is the Eilenberg-Moore construction for the comonad G. By a dual of Lemma 3.2, there is a pseudomonoid structure on A such that the whole triangle lifts to the Eilenberg-Moore construction in the bicategory MonMat(Set).

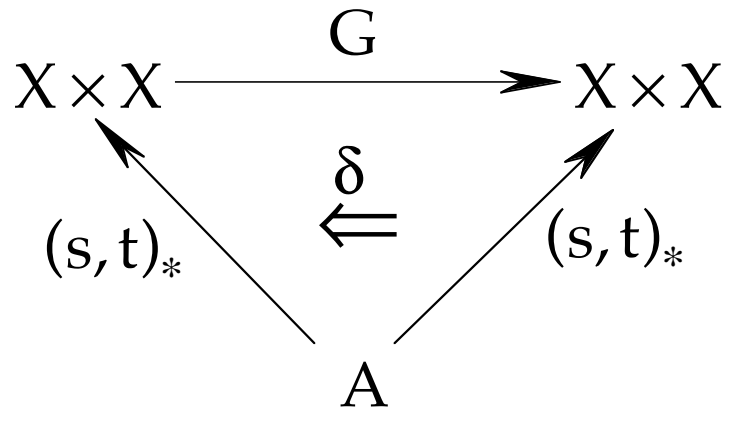

We already pointed out in the Introduction what the pseudomonoidal structure on A is; that on $X \times X$ is the special case of a chaotic category. Referring to the definition of basic data at the end of Section 9, we have:

**Proposition 10.1** *An equivalent definition of ordinary small categories is that they are basic data in the autonomous monoidal bicategory* $\mathrm{Mat(Set)}^{co}$ *where the special morphisms are all the maps.*

**Proof** Reversing 2-cells interchanges left and right adjunctions. So for a morphism to have a right adjoint in $\mathrm{Mat(Set)}^{co}$ is to be a left adjoint in Mat(Set); that is, to be the reverse of a function. Basic data in $\mathrm{Mat(Set)}^{co}$ therefore consists of a set X, a pseudomonoid A in Mat(Set), and a function $(s,t) : A \longrightarrow X \times X$ that is strong monoidal. The functor $\mathcal{B}(I, h^*)$ as at the end of Section 9 transports to the left-adjoint functor

$$\Sigma_{(s,t)} : Set / A \longrightarrow Set / X \times X$$

of the Introduction, which by Section 9 is strong monoidal. So we have a category **A**. Conversely, if **A** is a category, clearly $(s, t)$ is strong monoidal. **QED**

The discussion of the Introduction already shows that, if **A** is a groupoid, then it is $*$-autonomous in Mat(Set) with $Sa = a^{-1}$. In particular (the chaotic case), the endohom $X \times X$ is $*$-autonomous with $S(x,y) = (y,x)$. For **A** a groupoid, $(s,t)^* : X \times X \longrightarrow A$ is a strong $*$-autonomous map in $\mathrm{Mat(Set)}^{co}$. So we have Hopf basic data in $\mathrm{Mat(Set)}^{co}$. The converse almost holds.

**Proposition 10.2** *Consider a category as basic data in* $\mathrm{Mat(Set)}^{co}$. *The category is a groupoid iff the basic data are Hopf.*

**Proof** The characterizing property of $S = S_\ell$ is that

$$b \circ a = Sc \quad \text{iff} \quad c \circ b = S^{-1}a .$$

For each object $x$, put $e_x = S1_x$. Taking $c = 1_x$ and $b = S^{-1}a$ to ensure $c \circ b = S^{-1}a$, we deduce that $S^{-1}a \circ a = e_x$ for all $a : x \longrightarrow y$. Taking $a = e_x$ we see that $e_x = S^{-1}1_x$ so $S e_x = 1_x$. Now go back to the characterizing property with $c = e_x$, $b$ arbitrary, and $a = S(e_x \circ b)$ to ensure $c \circ b = S^{-1}a$: so we deduce that $b \circ S(e_x \circ b) = S e_x = 1_x$. It follows that every morphism $b$ has a right inverse. So the category is a groupoid. **QED**

**Remark 10.3** (This arose in lunchtime conversation with John Baez and Isar Stubbe.) The operation $S_\ell$ of $*$-autonomy is not unique. For a groupoid **A** as we have been considering, we can choose any endomorphism $e_x$ of each object $x$ and define $S_\ell a = a^{-1} \circ e_x$ so that $S_r a = e_x \circ a^{-1}$. This defines another $*$-autonomous structure on our pseudomonoid A.

## 11. Hopf bialgebroids

A bialgebroid A based on a k-algebra R is an opmonoidal monad on $R^e$ in $\mathrm{Mod}(\mathcal{V})$ (see Section 3). We have already seen that A becomes a k-algebra and that $\eta^* : A \longrightarrow R^e$ provides the Eilenberg-Moore object for the monad, thereby lifting to the bicategory of pseudomonoids in $\mathrm{Mod}(\mathcal{V})$.

In the terminology of Section 9, a bialgebroid is precisely basic data $\eta : R^e \longrightarrow A$ in $\mathcal{B} = \mathrm{Mod}(\mathcal{V})$. We define a bialgebroid $\eta : R^e \longrightarrow A$ to be *Hopf* when this basic data in $\mathrm{Mod}(\mathcal{V})$ is Hopf; that is, A should be $*$-autonomous and $\eta^* : A \longrightarrow R^e$ should be strong $*$-autonomous. It follows from Section 9 that $\mathrm{Mod}(\mathcal{V})(k, \eta^*)$ is strong monoidal and strong closed; this is none other than the functor

$$\mathcal{V}^A \longrightarrow \mathcal{V}^{R^e}$$

defined by restriction along $\eta : R^e \longrightarrow A$; compare Proposition 7.5 in the case of one-object $\mathcal{V}$-categories.

Preservation of internal homs was taken as paramount in the Hopf algebroid notions of [DS1] and [Sch]. Example 7.4 explains the connection between our work here and that of [DS1] while we see from the last paragraph that our Hopf bialgebroids are more restrictive than the Hopf algebroids of [Sch].

## 12. Quantum categories and quantum groupoids

It remains to state the main definitions of the paper. We now have the motivation and concepts readily at hand.

Let $\mathcal{V}$ be a braided monoidal category. We begin by recalling the definition of the right autonomous monoidal bicategory $\mathrm{Comod}(\mathcal{V})$ as appearing in the [DMS]. As there, we assume the condition:

*each of the functors* $X \otimes - : \mathcal{V} \longrightarrow \mathcal{V}$ *preserves equalizers.*

Briefly, $\mathrm{Comod}(\mathcal{V}) = \mathrm{Mod}(\mathcal{V}^{\mathrm{op}})^{\mathrm{coop}}$. To make calculations we will need to make the definition more explicit.

The objects of $\mathrm{Comod}(\mathcal{V})$ are comonoids C in $\mathcal{V}$; the comultiplication and counit are denoted by $\delta : C \longrightarrow C \otimes C$ and $\varepsilon : C \longrightarrow I$. The hom-category $\mathrm{Comod}(\mathcal{V})(C, D)$ is the category of Eilenberg-Moore coalgebras for the comonad $C \otimes - \otimes D$ on the category $\mathcal{V}$. This implies that the morphisms $M : C \longrightarrow\!\!\!\!\!\!\!+\;\; D$ in $\mathrm{Comod}(\mathcal{V})$ are comodules from C to D; that is, left C–, right D-comodules. So M is an object of $\mathcal{V}$ together with a coaction $\delta : M \longrightarrow C \otimes M \otimes D$ satisfying the expected equations. It is sometimes useful to deal with the left and right actions $\delta_\ell : M \longrightarrow C \otimes M$ and $\delta_r : M \longrightarrow M \otimes D$ which are obtained from $\delta$ using the counit. The 2-cells $f : M \Rightarrow M' : C \longrightarrow\!\!\!\!\!\!\!+\;\; D$ in $\mathrm{Comod}(\mathcal{V})$ are morphisms $f : M \longrightarrow M'$ in $\mathcal{V}$ respecting the coactions.

Composition of comodules $M : C \longrightarrow\!\!\!\!\!\!\!+\;\; D$ and $N : D \longrightarrow\!\!\!\!\!\!\!+\;\; E$ is given by the equalizer

$$N \circ M = M \underset{D}{\otimes} N \longrightarrow M \otimes N \overset{\delta_r \otimes 1}{\underset{1 \otimes \delta_\ell}{\rightrightarrows}} M \otimes D \otimes N .$$

The identity comodule $1_C : C \longrightarrow\!\!\!\!\!\!\!+\;\; C$ is $C$ with the obvious coaction.

The remaining details describing $\mathrm{Comod}(\mathcal{V})$ as a bicategory should now be clear.

**Remarks** (a) When $\mathcal{V} = \mathrm{Set}$, it is readily checked that $\mathrm{Comod}(\mathcal{V})$ is biequivalent to Mat(Set).

(b) The main case that should be kept in mind is when $\mathcal{V}$ is the category of vector spaces over a field k; then the objects of $\mathrm{Comod}(\mathcal{V})$ are precisely k-coalgebras.

(c) If $\mathcal{V}$ itself is a $*$-autonomous monoidal category then the distinction between $\mathrm{Mod}(\mathcal{V})$ and $\mathrm{Comod}(\mathcal{V})$ evaporates.

(d) By the Chu construction, any complete cocomplete closed monoidal $\mathcal{V}$ can be embedded into a complete cocomplete $*$-autonomous monoidal $\mathcal{E} = \mathcal{V}^{\mathrm{op}} \otimes \mathcal{V}$ taking V to $(1, V)$ where 1 is the terminal object of $\mathcal{V}$. The embedding is strong monoidal and preserves colimits and connected limits. So we can take full advantage of remark (c) by working in $\mathcal{E}-\mathrm{Mod}$ and deducing results for both $\mathrm{Mod}(\mathcal{V})$ and $\mathrm{Comod}(\mathcal{V})$.

Returning to general $\mathcal{V}$, we note that each comonoid morphism $f : C \longrightarrow D$ determines a comodule $f_* : C \longrightarrow\!\!\!\!\!\!\!+\;\; D$ defined to be C together with the coaction

$$C \xrightarrow{\delta} C \otimes C \xrightarrow{\delta \otimes f} C \otimes C \otimes D,$$

and a comodule $f^* : D \longrightarrow\!\!\!\!\!\!\!+\;\; C$ defined to be C together with the coaction

$$C \xrightarrow{\delta} C \otimes C \xrightarrow{f \otimes \delta} D \otimes C \otimes C.$$

Notice that we have $\gamma_f : f_* \circ f^* \Rightarrow 1_D$ which is defined to be $f : C \longrightarrow D$ since $f_* \circ f^* = f^* \underset{C}{\otimes} f_* = C$ with coaction $C \xrightarrow{\delta} C \otimes C \xrightarrow{1 \otimes \delta} C \otimes C \otimes C \xrightarrow{f \otimes 1 \otimes f} D \otimes C \otimes D$. Also, $f_* \underset{D}{\otimes} f^* = f^* \circ f_*$ is the equalizer

$$f_* \underset{D}{\otimes} f^* \longrightarrow C \otimes C \overset{(C \otimes f \otimes C) \circ (\delta \otimes C)}{\underset{(C \otimes f \otimes C) \circ (C \otimes \delta)}{\rightrightarrows}} C \otimes D \otimes C \quad ;$$

and, since

$$C \longrightarrow C \otimes C \overset{\delta \otimes C}{\underset{C \otimes \delta}{\rightrightarrows}} C \otimes C \otimes C$$

is an (absolute) equalizer, we have a unique morphism $C \longrightarrow f_* \underset{D}{\otimes} f^*$ commuting with the morphisms into $C \otimes C$; this gives us $\omega_f : 1_C \Rightarrow f^* \circ f_*$ . Indeed, $\gamma_f$ , $\omega_f$ are the counit and unit for an adjunction $f_* \dashv f^*$ in the bicategory $\mathrm{Comod}(\mathcal{V})$.

The comodules $f^*$ provide the special maps for the bicategory $\mathrm{Comod}(\mathcal{V})^{co}$.

Suppose C, D are comonoids. Then $C \otimes D$ becomes a comonoid with coaction

$$C \otimes D \xrightarrow{\delta \otimes \delta} C \otimes C \otimes D \otimes D \xrightarrow{C \otimes c_{C,D} \otimes D} C \otimes D \otimes C \otimes D$$

where c is the braiding and, as justified by coherence theorems, we ignore associativity in $\mathcal{V}$. For comodules $M : C \longrightarrow\!\!\!\!\!\!\!+\;\; C'$ and $N : D \longrightarrow\!\!\!\!\!\!\!+\;\; D'$, we obtain a comodule $M \otimes N : C \otimes D \longrightarrow\!\!\!\!\!\!\!+\;\; C' \otimes D'$ where the coaction is given in the obvious way using the braiding. This extends to a pseudofunctor $\otimes : \mathrm{Comod}(\mathcal{V}) \times \mathrm{Comod}(\mathcal{V}) \longrightarrow \mathrm{Comod}(\mathcal{V})$. The remaining structure required to obtain $\mathrm{Comod}(\mathcal{V})$ as a monoidal bicategory should be obvious.

Write $C^\circ$ for C with the comultiplication obtained from that of C by composing with the braiding. There is a pseudonatural equivalence between the category of comodules $M : C \otimes D \longrightarrow\!\!\!\!\!\!\!+\;\; E$ and the category of comodules $\hat{M} : D \longrightarrow\!\!\!\!\!\!\!+\;\; C^\circ \otimes E$, where $M = \hat{M}$ as objects. It follows that $C^\circ$ is a right bidual for C. This gives the structure of right autonomous monoidal bicategory to $\mathrm{Comod}(\mathcal{V})$.

Each $C^\circ \otimes C$ has the canonical structure of a pseudomonoid in $\mathrm{Comod}(\mathcal{V})$ because it is an endohom in the autonomous monoidal bicategory.

A *quantum category* over $\mathcal{V}$ is basic data C, $h : C^\circ \otimes C \longrightarrow A$ in $\mathrm{Comod}(\mathcal{V})^{co}$.

A *quantum groupoid* over $\mathcal{V}$ is Hopf basic data in $\mathrm{Comod}(\mathcal{V})^{co}$.

Centre of Australian Category Theory
Macquarie University
N. S. W. 2109
AUSTRALIA
email: street@math.mq.edu.au